\documentclass[reqno]{amsart}  
\theoremstyle{plain}
\newtheorem{theorem}{Theorem}[section]
\newtheorem{lemma}[theorem]{Lemma}
\newtheorem{remark}[theorem]{Remark}
\newtheorem{proposition}[theorem]{Proposition}
\newtheorem{problem}[theorem]{Problem}

\newtheorem{corollary}[theorem]{Corollary}
\newtheorem{example}[theorem]{Example}

\numberwithin{equation}{section}
\allowdisplaybreaks[1]

\theoremstyle{definition}

\newtheorem{definition}[theorem]{Definition}

\theoremstyle{remark}

\newcommand{\bA}{{\mathbf A}}
\newcommand{\bU}{{\mathbf U}}

\newcommand{\blam}{{\boldsymbol \lambda}}
\newcommand{\lam}{\lambda}
\newcommand{\bzeta}{{\boldsymbol \zeta}}

\newcommand{\bn}{{\mathbf n}}

\newcommand{\cD}{{\mathcal D}}
\newcommand{\cF}{{\mathcal F}}
\newcommand{\cH}{{\mathcal H}}
\newcommand{\cL}{{\mathcal L}}

\newcommand{\cR}{{\mathcal R}}
\newcommand{\cS}{{\mathcal S}}

\newcommand{\cU}{{\mathcal U}}
\newcommand{\cX}{{\mathcal X}}
\newcommand{\cY}{{\mathcal Y}}
\newcommand{\cO}{{\mathcal O}}

\newcommand{\C}{{\mathbb C}}
\newcommand{\B}{{\mathbb B}}

\newcommand{\sbm}[1]{\left[\begin{smallmatrix} #1
		\end{smallmatrix}\right]}

\begin{document}

\title[Realization for multipliers]{Transfer-function realization for
multipliers of the Arveson space}
\author[J. A. Ball]{Joseph A. Ball}
\address{Department of Mathematics,
Virginia Tech,
Blacksburg, VA 24061-0123, USA}
\email{ball@math.vt.edu}
\author[V. Bolotnikov]{Vladimir Bolotnikov}
\address{Department of Mathematics,
The College of William and Mary,
Williamsburg VA 23187-8795, USA}
\email{vladi@math.wm.edu}
\author[Q. Fang]{Quanlei Fang}
\address{Department of Mathematics,
Virginia Tech,
Blacksburg, VA 24061-0123, USA}
\email{qlfang@math.vt.edu}

\begin{abstract}
An interesting and recently much studied generalization of the
classical Schur class is the class of contractive operator-valued
multipliers for the reproducing kernel Hilbert space $\cH(k_{d})$ on
the unit ball ${\mathbb B}^{d} \subset {\mathbb C}^{d}$, where
$k_{d}$ is the positive kernel $k_{d}(\blam, \bzeta) = 1/(1
- \langle \blam, \bzeta \rangle)$ on ${\mathbb B}^{d}$.  We study
this space from the point of view of realization theory and
functional models of de Branges-Rovnyak type.  We highlight features
which depart from the classical univariate case:  coisometric realizations
have only partial uniqueness properties, the nonuniqueness can be
described explicitly, and this description assumes a particularly concrete
form in the functional-model context.
\end{abstract}


\maketitle

\section{Introduction}  \label{S:Intro}
\setcounter{equation}{0}
Let $\cU$ and $\cY$ be two Hilbert spaces and
let $\cL(\cU, \cY)$ be the space of all bounded linear operators
between $\cU$ and $\cY$.  We also let $H^2_{\cU}$ be the standard Hardy
space of the $\cU$-valued holomorphic functions on the unit disk ${\mathbb
D}$. The operator-valued version of the classical Schur class ${\mathcal
S}(\cU, \cY)$ is defined to be the set of all holomorphic, contractive
$\cL(\cU, \cY)$-valued functions on ${\mathbb D}$.
The following equivalent characterizations of the Schur class are
well known.

\begin{theorem}  \label{T:clas-Schur}
           Let $ S \colon {\mathbb D} \to \cL(\cU, \cY)$ be given.
           Then the following are equivalent:
           \begin{enumerate}
\item[(1)] $S \in \cS(\cU, \cY)$, i.e., $S$ is holomorphic on
                 ${\mathbb D}$ with $\| S(\lam) \| \le 1$ for all $\lam \in
                 {\mathbb D}$.

\item[(1$^{\prime}$)]  The multiplication operator $M_{S} \colon
              f(z) \mapsto S(z) \cdot f(z)$ is a contraction from
              $H^{2}_{\cU}$ into $H^{2}_{\cY}$: $\|M_{S}\|_{\text{op}} \le 1$.

	\item[(2)] The associated kernel function
	\begin{equation}  \label{clas-KS}
	K_{S}(\lam, \zeta) = \frac{ I_{\cY} - S(\lam) S(\zeta)^{*}}{1 -
	\lam \overline{\zeta}}
	\end{equation}
	is a positive kernel on ${\mathbb D}\times {\mathbb D}$, i.e.,
              there exists an	operator-valued function $H \colon
              {\mathbb D} \to \cL(\cH,\cY)$ for some auxiliary Hilbert space
              $\cH$ so that
	\begin{equation}  \label{clas-KSfact}
	  K_{S}(\lam, \zeta) = H(\lam) H(\zeta)^{*}.
	\end{equation}

	\item[(3)] There is an auxiliary Hilbert space $\cX$ and a
         unitary connecting operator
         $$
         U = \begin{bmatrix} A & B \\ C & D \end{bmatrix} \colon
         \begin{bmatrix} \cX \\ \cU \end{bmatrix} \to \begin{bmatrix} \cX \\
             \cY \end{bmatrix}
         $$
         so that $S(\lam)$ can be expressed as
         \begin{equation}  \label{clas-real}
           S(\lam) = D + \lam C (I - \lam A)^{-1} B.
           \end{equation}

           \item[(4)] $S(\lam)$ has a realization as in \eqref{clas-real} where
           the connecting operator $U$ is any one of (i) isometric, (ii)
	coisometric, or (iii) contractive.
        \end{enumerate}
        \end{theorem}

        We remark that the proof that the coisometric version of (4) implies
        (2) in Theorem \ref{T:clas-Schur} is particularly transparent: if
        $S(\lam)$ has the form \eqref{clas-real} with $U = \sbm{ A & B  \\C &
        D }$ coisometric, a simple calculation reveals that
        \eqref{clas-KSfact} holds with $H(\lam) = C ( I - \lam A)^{-1}$,
        i.e.,
        \begin{equation}  \label{clas-KSreal}
            K_{S}(\lam, \zeta) = C (I - \lam A)^{-1} (I - \overline{\zeta}
            A^{*})^{-1} C^{*}:= K_{C,A}(\lam, \zeta).
        \end{equation}

        Among all the possible classes for the connecting operator $U$ (i.e.,
        unitary, isometric, coisometric or simply contractive), the class of
        coisometric ones is particularly prominent due to its connection
        with functional-model realizations using the de Branges-Rovnyak
        reproducing kernel Hilbert space $\cH(K_{S})$
        associated with the positive kernel $K_{S}$ given by \eqref{clas-KS}.
        We recall (see the original work of Aronszajn \cite{aron}) that any
        positive kernel $(\lam, \zeta) \mapsto k(\lam, \zeta) \in \cL(\cY)$
        on a set $\Omega \times \Omega$ (so $\lam, \zeta \in \Omega$)
        gives rise to  a reproducing kernel Hilbert space (RKHS) $\cH(k)$
        consisting of $\cY$-valued functions on $\Omega$ with the defining
        property: for each $\zeta \in \Omega$ and $y \in \cY$, the
            $\cY$-valued function $(k_{\zeta} y)(\lambda):=k(\lam, \zeta) y$
           is in $\cH(k)$ and has the reproducing property
           $$
            \langle f, k_{\zeta}y \rangle_{\cH(k)} = \langle f(\zeta), y
            \rangle_{\cY}\quad \text{for all}\quad y \in \cY, \; f\in\cH(k).
           $$
          We remark that the Hardy space $H^{2}_{\cY}$ is the RKHS
          associated with the Szeg\"o kernel $k_{\text{Sz}}(\lam,
          \zeta) =(1 - \lam \overline{\zeta})^{-1}I_{\cY}$ positive on
          ${\mathbb D} \times {\mathbb D}$ where ${\mathbb D}$ is the unit
          disk.   Applying Aronszajn's construction
         to the positive kernel $K_{S}$ on ${\mathbb D}$ for a
         Schur-class function $S$ as in \eqref{clas-KSreal} gives the
         reproducing kernel Hilbert space $\cH(K_{S})$, the de
         Branges-Rovnyak space associated with $S$.
         Then we have the following concrete, functional-model realization for
         the Schur-class function $S$ \cite{dbr1, dbr2}.

         \begin{theorem}  \label{T:clas-dBRreal}
            Suppose that $S \in {\mathcal S}(\cU, \cY)$ and let $\cH(K_{S})$
            be the associated de Branges-Rovnyak model space.  Then
            the connecting operator
            $$
            U = \begin{bmatrix} A & B \\ C & D \end{bmatrix} \colon
            \begin{bmatrix} \cH(K_{S}) \\ \cU \end{bmatrix} \to
	 \begin{bmatrix} \cH(K_{S}) \\ \cY \end{bmatrix}
             $$
           with the entries defined by
\begin{align}
               A \colon f(\lam) \mapsto \frac{f(\lam) - f(0)}{\lam}, &
	     \qquad \quad
C \colon f \mapsto f(0) \quad\mbox{for} \; \; f
\in\cH(K_{S}),\label{clas-dBRrealac} \\
B \colon u \mapsto \frac{S(\lam) - S(0)}{\lam}u,  & \qquad
\quad D \colon u \mapsto S(0) u \quad\mbox{for} \; \; u\in \cU,
\label{clas-dBRreal}
           \end{align}
           provides a coisometric realization of $S(\lam)$, i.e.,
           $U$ is coisometric as an operator from $\sbm{\cH(K_{S}) \\ \cU }$
           to $\sbm{\cH(K_{S}) \\ \cY}$ and we recover $S(\lam)$ via the formula
           \eqref{clas-real}.
           \end{theorem}
The de Branges-Rovnyak functional-model realization is {\em closely
outer-connected} in the sense that the pair $(C,A)$ is {\em observable},
i.e., that
$$
C(I - zA)^{-1} x = 0 \quad\text{for all} \; \; z \in {\mathbb D}
\quad\Longrightarrow \quad x = 0.
$$
Observability of the pair $(C,A)$ is a minimality condition under which
the coisometric realization is essentially unique: {\em every coisometric
closely outer-connected realization of an $S\in\cS(\cU,\cY)$ is unitarily
equivalent to the de Branges-Rovnyak functional-model realization}.
It can also be shown that, if $(C,A)$ is observable and if $U = \sbm{
A & B \\ C &
S(0) }$ provides a coisometric realization for the $S(\lam)$, then the
operator $B$ is already uniquely determined by $C,A$ and $S$ (see
Remark \ref{R:collapse} below).

A multivariable generalization of the Szeg\"o kernel much studied of
late is the positive kernel
$$
k_d(\blam,\bzeta)=\frac{1}{1-\langle \blam,  \bzeta  \rangle}
$$
on ${\mathbb B}^{d} \times {\mathbb B}^{d}$ where
$\B^d=\left\{\blam=(\lam_1,\dots, \lam_d)\in\C^d \colon  \langle \blam,
\blam\rangle<1\right\}$ is the unit ball of the  $d$-dimensional Euclidean
space $\C^d$. By
$$
\langle \blam,  \bzeta \rangle=\langle \blam,
\bzeta  \rangle_{\C^d}=\sum_{j=1}^d \lam_j \overline{\zeta}_j\quad
\text{for} \; \; \blam, \bzeta\in\C^d
$$
we mean the standard inner product in $\C^d$.
The associated RKHS $\cH(k_d)$ obtained via Aronszajn's
construction is a natural multivariable analogue
of the Hardy space $H^2$ of the unit disk and coincides with $H^2$ if
$d=1$.

For $\cY$ an auxiliary Hilbert space, we consider the tensor product
Hilbert space $\cH_\cY(k_d):=\cH(k_d)\otimes\cY$ whose
elements can be viewed as  $\cY$-valued functions in $\cH(k_d)$.
The space of multipliers ${\mathcal M}_d(\cU,\cY)$ is defined  as the
space of all
$\cL(\cU,\cY)$-valued analytic functions  $S$ on $\B^d$ such that the
induced multiplication operator
\begin{equation}
M_S: \; f(\blam)\to S(\blam)\cdot f(\blam)
\label{ms}
\end{equation}
maps $\cH_{\cU}(k_d)$ into $\cH_{\cY}(k_d)$. It follows by the closed
graph theorem that for every $S\in {\mathcal M}_d(\cU,\cY)$, the
operator $M_S$ is  bounded. We shall pay particular attention to the unit
ball of  $ {\mathcal M}_{d}(\cU, \cY)$, denoted by
$$
{\mathcal S}_{d}(\cU, \cY) = \{ S \in {\mathcal M}_{d}(\cU,\cY) \colon
\| M_{S} \|_{\text{op}} \le 1 \}.
$$
Since ${\mathcal S}_{1}(\cU, \cY)$ collapses to
the classical Schur class (by the equivalence
(1) $\Longleftrightarrow$ (1$^{\prime}$) in Theorem
\ref{T:clas-Schur}),  we refer to ${\mathcal S}_{d}(\cU,
\cY)$ as a generalized ($d$-variable) {\em Schur class}.
The following result appears in \cite{BTV,
AM} and is the precise
analogue of Theorem \ref{T:clas-Schur} for the multivariable case.
Note that there is no analogue of condition (1) in Theorem
\ref{T:clas-Schur} and condition (1) in Theorem \ref{T:BTV} is the
analogue of condition (1$^{\prime}$) in Theorem \ref{T:clas-Schur}.

\begin{theorem}
\label{T:BTV}
Let $S$ be an $\cL(\cU, \cY)$-valued function defined
on ${\mathbb B}^{d}$.  The following are equivalent:
\begin{enumerate}
\item $S$ belongs to $\cS_d(\cU, \cY)$.
\item The kernel
\begin{equation}  \label{KS}
         K_{S}(\blam, \bzeta) = \frac{I_{\cY} - S(\blam) S(\bzeta)^{*}}{1 -
         \langle \blam, \bzeta \rangle}
\end{equation}
is positive on $\B^{d}\times\B^{d}$.
\item There exists a Hilbert space $\cX$ and a
unitary connecting operator (or colligation) $\bU$ of the form
\begin{equation}
\label{1.7a}
\bU = \begin{bmatrix} A  & B \\ C & D \end{bmatrix} =
\begin{bmatrix} A_{1} & B_{1} \\ \vdots & \vdots \\ A_{d} & B_{d}
\\ C & D \end{bmatrix} \colon \begin{bmatrix} \cX \\ \cU
\end{bmatrix} \to \begin{bmatrix}\cX^{d} \\ \cY \end{bmatrix}
\end{equation}
so that $S(\blam)$ can be realized in the form
\begin{eqnarray}
S(\blam)&=&D+C\left(I_{\cX}-\lam_1A_1-\cdots-\lam_dA_d\right)^{-1}
\left(\lambda_1B_1+\ldots+\lambda_dB_d\right)\nonumber\\
&=& D + C (I - Z(\blam) A)^{-1} Z(\blam) B
\label{1.5a}
\end{eqnarray}
where we set
\begin{equation}
Z(\blam)=\begin{bmatrix}\lam_1 I_{\cX} & \ldots &
\lam_dI_{\cX}\end{bmatrix},\quad
A=\begin{bmatrix} A_1 \\ \vdots \\ A_d\end{bmatrix},\quad
B=\begin{bmatrix} B_1 \\ \vdots \\ B_d\end{bmatrix}.
\label{1.6a}
\end{equation}
\item There exists a Hilbert space $\cX$ and a
contractive connecting operator $\bU$ of the form \eqref{1.7a}
so that $S(\blam)$ can be realized in the form \eqref{1.5a}.
\end{enumerate}
\end{theorem}
Although Statement (4) in Theorem \ref{T:BTV} concerning contractive
realizations does not appear in \cite{AM, BTV}, its equivalence to statements
(1)-(3) is quite obvious. Indeed, implication $(3)\Rightarrow(4)$ is
trivial; on the other hand, a straightforward calculation (see e.g.,
\cite[Lemma 2.2]{ADR}) shows that for $S$ of the form \eqref{1.5a},
\begin{align}
        K_S(\blam, \bzeta)&=C (I_{\cX} - Z(\blam) A)^{-1}(I_{\cX} -
A^{*}Z(\bzeta)^{*})^{-1} C^{*}\label{1.6u}\\
&+\begin{bmatrix}C (I-Z(\blam)A)^{-1}Z(\blam) & I\end{bmatrix}
\frac{I-\bU\bU^*}{1-\langle \blam, \bzeta \rangle}
\begin{bmatrix}Z(\bzeta)^*(I_{\cX} -
A^{*}Z(\bzeta)^{*})^{-1} C^{*}\\ I\end{bmatrix}
\notag
\end{align}
where $\bU$ is defined in \eqref{1.7a}. Thus, if $\bU$ is a
contraction,
the kernel $K_S(\blam, \bzeta)$ is positive on $\B^{d}\times\B^{d}$
which proves implication $(4)\Rightarrow(2)$ in Theorem \ref{T:BTV}.

In analogy with the univariate case, a realization of the form
(\ref{1.5a}) is called {\em coisometric}, {\em isometric}, {\em unitary}
or {\em contractive} if the operator $\bU$ is respectively, coisometric,
isometric, unitary or just contractive. It turns out that a more useful
analogue of ``coisometric realization'' appearing in the classical univariate
case is not that the whole connecting operator ${\bf U}^{*}$ be isometric,
but rather that ${\bf U}^{*}$ be isometric on a certain canonical subspace
of $\cX^d\oplus\cY$.
\begin{definition}
A realization \eqref{1.5a} of $S\in\cS_d(\cU,\cY)$ is called
{\em weakly coisometric} if the adjoint $\bU^*: \, \cX^d\oplus\cY\to
\cX\oplus\cU$ of the connecting operator
is contractive and isometric on the subspace
$\begin{bmatrix}{\mathcal D}_{C,A}\\
\cY\end{bmatrix}\subset \begin{bmatrix}\cX^d\\ \cY\end{bmatrix}$ where
\begin{equation}  \label{domV0}
         \cD=\cD_{C,A} := \operatorname{\overline{span}}
         \{Z(\bzeta)^{*}(I - A^{*} Z(\bzeta)^{*})^{-1} C^{*}y \colon \; \;
\bzeta\in {\mathbb B}^{d}, \, y \in \cY \} \subset \cX^{d}.
         \end{equation}
\label{D:1.2}
\end{definition}
The notion of weakly coisometric realizations has been introduced in
\cite{BTV}. It does not appear in the single-variable context for a
simple reason that if the pair  $(C,A)$ is observable, then
a weakly coisometric realization is automatically
coisometric (see \cite[p.~100]{BTV} and also Remark \ref{R:collapse}
below). The following intrinsic kernel
characterization as to when a given contractive realization is a
weakly coisometric realization turns out to be a convenient tool for
our current purposes. Equality \eqref{2.1} below is the multivariable
analogue of equality \eqref{clas-KSreal}.
\begin{proposition} \label{R:1.3}
A contractive realization \eqref{1.5a} of $S\in\cS_d(\cU,\cY)$ is  weakly
coisometric if and only if the kernel  $K_{S}(\blam, \bzeta)$ associated
to $S$ via \eqref{KS}  can alternatively be written as
\begin{equation}
K_{S}(\blam, \bzeta) = K_{C, \bA}(\blam, \bzeta)
\label{2.1}
\end{equation}
where
\begin{equation}  \label{KCA}
         K_{C,{\mathbf A}}(\blam, \bzeta):=
         C (I - \lam_{1}A_{1} - \cdots -   \lam_{d} A_{d})^{-1}
         (I - \overline{\zeta_{1}} A_{1}^{*} - \cdots - \overline{\zeta_{d}}
         A_{d}^{*})^{-1} C^{*}.
\end{equation}
\end{proposition}
\begin{proof} Let $\bU = \sbm{ A & B \\ C & D}$ be the connecting
operator of a  contractive realization of $S\in\cS_d(\cU,\cY)$.
It is readily seen from the formula \eqref{1.6u}, that equality
\eqref{2.1} holds if and only if the operator $\bU^*$ is isometric on the
space
$$
{\mathcal M}:=\overline{\operatorname{span}}\left\{\begin{bmatrix}
Z(\bzeta)^*(I_{\cX}-A^{*}Z(\bzeta)^{*})^{-1} C^{*}\\ I\end{bmatrix}y
\colon \; \;  \bzeta \in\B^d, \; y\in\cY \right\}\subset \cX^d\oplus \cY.
$$
By setting $\bzeta=0$ in the last formula, we see that $\sbm{ 0 \\ y }
\in{\mathcal M}$ for  all $y \in \cY$ and thus ${\mathcal M}$ splits in
the form ${\mathcal M}=\sbm{{\mathcal D}\\ \cY}$ where $\cD$ is defined in
\eqref{domV0}. The rest follows by Definition \ref{D:1.2}.\end{proof}

The present paper analyzes a number of finer structural issues
surrounding a Schur-class function $S(\blam)$ and its associated
positive kernel \eqref{KS}.  We analyze when equality \eqref{2.1}
holds in both a realization and a purely function-theoretic context.
We analyze the problem of realizing a kernel of the form
$K_{C,{\mathbf A}}(\blam, \bzeta)$  as $K_{S}(\blam, \bzeta)$ for a
Schur-class function $S \in {\mathcal S}_d(\cU, \cY)$ (Theorems
\ref{T:CAtoSa} and \ref{T:repr}) and we analyze the
nonuniqueness of the input operator $B$ inherent in a weakly coisometric
(as well as coisometric or unitary)
realization of a given Schur-class function $S \in {\mathcal S}_{d}(\cU, \cY)$
using a given output-pair $(C,{\bf A})$ which is observable in an
appropriate multivariable sense (Theorem \ref{T:CAtoS} and Theorem
\ref{T:3.5}).
Upon applying Aronszajn's construction to the kernel $K_S$ associated
with a Schur-class function $S\in\cS_d(\cU,\cY)$ (which is positive on
$\B^d$ by Theorem \ref{T:BTV}), one gets the  de Branges-Rovnyak space
$\cH(K_{S})$ that can serve as the state space for a weakly coisometric
realization for $S$. A weakly coisometric realization for $S$ with the
state space equal to $\cH(K_S)$ and with the output operator $C$ equal to
evaluation at zero on $\cH(K_{S})$ will be called a {\em generalized
functional-model realization}\footnote{The term (not necessarily
generalized) {\em functional-model realization} is explained below.}.

Our earlier paper \cite{BBF1} focuses on the structure of reproducing
kernel Hilbert spaces $\cH(K_{C, \bA})$ with reproducing kernel
$K_{C, \bA}$ of the form \eqref{KCA}.  Such spaces can be viewed as
the range of an observability operator associated with a state-output
multidimensional linear system of the form
$$ \Sigma \colon \left\{ \begin{array}{rcl}
x(\bn) & = & A_{1} x(\sigma_{1}(\bn)) + \cdots + A_{d}
        x(\sigma_{d}(\bn)) \\
       y(\bn) & =  & C x(\bn)
\end{array} \right.
$$
where
$$
\sigma_{k}(\bn) = \sigma_{k}((n_{1}, \dots, n_{d})) =
(n_{1}, \dots, n_{k-1}, n_{k}+1, n_{k+1}, \dots, n_{d})
$$
for $\bn = (n_{1}, \dots, n_{d}) \in {\mathbb Z}^{d}_{+}.$
Also discussed in \cite{BBF1} are connections with noncommutative
analogues of these objects, where the reproducing kernel Hilbert space
is of the noncommutative type discussed in \cite{NFRKHS} consisting
of formal power series with vector coefficients and where the system
has evolution along a free semigroup rather than along ${\mathbb
Z}^{d}_{+}$.
The paper \cite{BBF1} also serves as a resource for the present paper,
since, once one has established the equality \eqref{2.1}, results concerning
$\cH(K_{C, \bA})$ from \cite{BBF1} immediately
yield the corresponding result for the space $\cH(K_{S})$.

We reserve the term (non-generalized)
{\em functional-model realization} for the case where $\cH(K_{S})$ is
invariant under the adjoints $M_{\lambda_{j}}^{*}$ of the
multiplication operators $M_{\lambda_{j}} \colon f(\blam) \mapsto
\lambda_{j} f(\blam)$ on $\cH_{\cY}(k_{d})$ and the state-space
operators $\bA = (A_{1}, \dots, A_{d})$ in the realization are taken
to be $A_{j} = M_{\lambda_{j}}^{*}|_{\cH(K_{S})}$;  the
characteristic function $S_{{\mathbf T}}(\blam)$ for a commuting row
contraction ${\mathbf T} = (T_{1}, \dots, T_{d})$ (see \cite{BES,
BES2, BS}) as well as inner
functions (Schur-class multipliers $S$ for which the associated
multiplication operator
$M_{S} \colon f(\blam) \mapsto S(\blam) \cdot f(\blam)$ is a partial
isometry) are of this type.  We discuss the special features of this
case (where $\cH(K_{S})$ is invariant under $M_{\lambda_{j}}^{*}$ for
$j = 1, \dots, d$ and where $S(\blam)$ has a realization with
commuting state-space operators $A_{1}, \dots, A_{d}$)
in our separate paper \cite{BBF2b}.

The paper is organized as follows.
Section 2 develops the ideas surrounding observable weakly coisometric
realizations and the quantification of the nonuniqueness of the
input operator in such realizations.
In Section 3 we show that any
Schur-class function $S\in\cS_d(\cU,\cY)$ admits a generalized functional-model
realization and that any observable  weakly coisometric realization of
$S$ is unitarily equivalent to some generalized functional-model
realization. Preliminary results of this latter type appear in the
paper of Alpay-Dijksma-Rovnyak \cite{ADR}.  In Section 4 we introduce
a general setting for the overlapping spaces appearing prominently in
the work of de Branges and Rovnyak \cite{dbr1, dbr2} and indicate how
special cases of these spaces appear in Sections 2 and 3 in
connection with the nonuniqueness of the input
operator in observable weakly coisometric realizations.

In our followup paper \cite{BBF3}, we develop the noncommutative
theory parallel to the results of the present paper.  In this
setting, the Schur-class function $S$ becomes a formal power series
in noncommuting indeterminates inducing a contractive multiplication
operator between Fock-Hilbert spaces consisting of formal power
series with vector coefficients.  Such a Schur-class multiplier
induces a kernel $K_{S}(z,w)$ in noncommuting indeterminates $z =
(z_{1}, \dots, z_{d})$ and $w = (w_{1}, \dots, w_{d})$ which is a
noncommutative positive kernel in the sense of \cite{NFRKHS}.  The
associated noncommutative formal reproducing kernel Hilbert space
$\cH(K_{S})$ is a noncommutative analogue of the space $\cH(K_{S})$
studied here (where elements of the space of functions of commuting
variables $\blam = (\lam_{1}, \dots, \lam_{d})$ and is an alternative
multivariable generalization of the classical case
(\cite{dbr1, dbr2}).  For this setting the analogy with the classical
case turns out to be more compelling than for the case of several
commuting variables presented here.

\section{Weakly coisometric realizations}
\label{S:coisom}

Weakly coisometric realizations of Schur class functions are closely
related to range spaces of observability operators studied in \cite{BBF1}.
Let $\bA=(A_1,\ldots,A_d)$ be a $d$-tuple of operators in $\cL(\cX)$.
If $C \in \cL(\cX, \cY)$, then the pair $(C,{\mathbf A})$ is said
         to be an {\em output pair}.  Such an output pair is said to be {\em
         contractive} if
         $$
A_{1}^{*}A_{1} + \cdots + A_{d}^{*} A_{d} + C^{*}C\le I_{\cX},
         $$
to be {\em isometric} if equality holds in the above relation, and
to be {\em output-stable} if the associated observability
         operator
         \begin{equation}  \label{cO}
           \cO_{C, \bA} \colon x \mapsto C (I - \lam_{1} A_{1} - \cdots -
           \lam_{d} A_{d})^{-1} x
          \end{equation}
maps $\cX$ into $\cH_{\cY}(k_{d})$. As it was shown in \cite{BBF1}, any
contractive pair $(C,\bA)$ is output stable and, moreover, the corresponding
observability operator ${\mathcal O}_{C, \bA}: \, \cX\to\cH_{\cY}(k_d)$ is a
contraction. An output stable pair $(C,\bA)$ is called {\em observable} if
the observability operator
${\mathcal O}_{C, \bA}$ is injective, i.e.,
$$
C (I- \lam_{1} A_{1} - \cdots -
\lam_{d} A_{d})^{-1}x\equiv 0\quad\Longrightarrow \quad
x=0.
$$
The following result from \cite{BBF1} gives the close connection
between spaces of the form $\cH(K_{C, \bA})$ and ranges of
observability operators.

\begin{theorem}  \label{T:3-1.2nc}
Let $(C,\bA)$ be a contractive pair with $C\in{\mathcal L}(\cX,\cY)$
and with associated positive kernel $K_{C, \bA}$ given by \eqref{KCA} and
the observability operator $\cO_{C, \bA}$ given by \eqref{cO}.  Then:
         \begin{enumerate}
\item The reproducing kernel Hilbert space $\cH(K_{C,\bA})$ is
characterized as
$$
\cH(K_{C,\bA}) = \operatorname{Ran}\, \cO_{C, \bA}
$$
with the lifted norm given by
$\; \| \cO_{C, \bA} x \|_{\cH(K_{C, \bA})} = \| Q x \|_{\cX}$,
where $Q$ is the orthogonal projection onto
$(\operatorname{Ker}\, \cO_{C, \bA})^{\perp}$.
        \item The operator ${\cO}_{C,\bA}$
               is a contraction of $\cX$ into $\cH(K_{C,\bA})$.
               It is an isometry if and only if the the pair $(C,\bA)$ is
               observable.

               \item There exist operators $T_1,\ldots,T_d\in{\mathcal
                L}(\cH(K_{C,\bA}))$ such that relations
$$
               f(\blam)-f(0)=\sum_{j=1}^d \lambda_j(T_jf)(\blam)\quad(\blam
               \in\B^d)
$$
               and
$$             \sum_{j=1}^d \|T_jf\|^2_{\cH(K_{C,\bA})}
               \le\|f\|^2_{\cH(K_{C,\bA})}-\|f(0)\|^2_{\cY}
$$
hold for every  function $f\in\cH(K_{C,\bA})$.
               \end{enumerate}
              \end{theorem}
Proposition \ref{R:1.3} and Theorem \ref{T:3-1.2nc} assert that every
weakly coisometric realization
of a Schur-class function $S\in\cS_d(\cU,\cY)$ identifies the corresponding de
Branges-Rovnyak space $\cH(K_S)$ as the range space of the
observability operator corresponding to a contractive pair $(C,\bA)$.
The next proposition shows that the reverse identification is also
possible.

        \begin{theorem}  \label{T:CAtoSa}
Let $(C, \bA)$ with $C\in\cL(\cX,\cY)$ be a contractive pair.  Then there
exist an input space $\cU$ and an $S\in\cS_{d}(\cU, \cY)$ such that
\begin{equation}
K_{S}(\blam, \bzeta) = K_{C, \bA}(\blam, \bzeta).
\label{2.1a}
\end{equation}
\end{theorem}

\begin{proof}  Choose a Hilbert space $\cU$ with
$$
\operatorname{dim}\, \cU \ge \operatorname{rank}\, \left(
\begin{bmatrix}I_{\cX^{d}} &
       0 \\  0 & I_{\cY}\end{bmatrix} - \begin{bmatrix}A \\ C\end{bmatrix}
       \begin{bmatrix} A^{*} & C^{*} \end{bmatrix}\right)
$$
       and let $\sbm{B \\ D} \colon \cU \to \cX^{d} \oplus \cY$
       be a solution of the Cholesky
       factorization problem
       $$
       \begin{bmatrix} B \\ D \end{bmatrix} \begin{bmatrix} B^{*} &
	D^{*}\end{bmatrix} = \begin{bmatrix} I_{\cX^{d}} & 0 \\ 0 &
	I_{\cY}
       \end{bmatrix} - \begin{bmatrix} A \\ C \end{bmatrix}
       \begin{bmatrix} A^{*} & C^{*} \end{bmatrix}.
      $$
      Then $\sbm{A & B  \\ C & D }$ is coisometric.  Let $S(\blam)$ be
      given by the realization formula \eqref{1.5a}. Then $S \in
      \cS_{d}(\cU, \cY)$ and Proposition
      \ref{R:1.3} guarantees \eqref{2.1a} as wanted.
      \end{proof}

Theorem \ref{T:CAtoSa} shows that every range space
$\operatorname{Ran}\, {\mathcal O}_{C, \bA}=\cH(K_{C,\bA})$ associated with a
contractive pair  $(C, \bA)$ can be considered as the de Branges-Rovnyak
space $\cH(K_S)$ for an appropriately chosen Schur-class function $S$, which we
will call {\em a representer of $\cH(K_{C,\bA})$}. A description of all
representers for a given $\cH(K_{C,\bA})$ will be given below in Theorem
\ref{T:repr}.

Now we discuss equality \eqref{2.1a} independently of the realization
context.
With a given contractive pair $(C, \bA)$ with $C\in\cL(\cX,\cY)$ and an
$\cL(\cU,\cY)$-valued function $S$ defined on $\B^d$  we associate the
operator
\begin{equation}
V = \begin{bmatrix} A_{V} & B_{V} \\ C_{V} & D_{V} \end{bmatrix}
           \colon \begin{bmatrix} \cD \\ \cY \end{bmatrix} \to
           \begin{bmatrix} \cX \\ \cU \end{bmatrix}
\label{defv}
\end{equation}
(where the space $\cD$ is defined in \eqref{domV0}) with the
entries given by
\begin{equation}
A_{V}=A^*\vert_{\cD},\qquad B_{V}=C^*,\qquad D_V=S(0)^*,
\label{abd}
\end{equation}
and where $C_{V}$ is uniquely determined by linearity and continuity
by its action on a generic generating vector for $\cD$:
\begin{equation}
C_{V} \colon \; Z(\bzeta)^{*} (I - A^{*}Z(\bzeta)^{*})^{-1}C^{*}y \
\mapsto (S(\bzeta)^{*} - S(0)^{*}) y\quad \mbox{for}\; \;
\bzeta\in\B^d\; \; y\in\cY.
\label{c}
\end{equation}

\begin{lemma}
Let  $(C, \bA)$ be a contractive pair and let $S$ be an
$\cL(\cU,\cY)$-valued function defined on $\B^d$. Then
\eqref{2.1a} holds (and therefore also $S$ belongs to $\cS_d(\cU,\cY)$) if and
only if the operator $V$ defined in
\eqref{defv}--\eqref{c} is an isometry from $\cD\oplus\cY$ onto
        \begin{equation}  \label{ranV}
           \cR_{V}:= \operatorname{\overline{span}} \left\{ \begin{bmatrix}
           (I - A^{*}Z(\bzeta)^{*})^{-1}C^{*}y \\ S(\bzeta)^{*}y
\end{bmatrix} \colon \bzeta \in
           {\mathbb B}^{d}, \, y \in \cY \right\} \subset \begin{bmatrix} \cX
           \\ \cU \end{bmatrix}.
        \end{equation}
\label{L:2.1}
\end{lemma}
\begin{proof} Let  $S\in\cS_d(\cU,\cY)$ and let equality \eqref{2.1a}
hold, i.e., let
$$
\frac{ I_{\cY} - S(\blam) S(\bzeta)^{*}}{1 - \langle \blam, \bzeta\rangle}
        = C (I - Z(\blam) A)^{-1}(I - A^{*} Z(\bzeta)^{*})^{-1}C^{*},
$$
which can be written equivalently (due to the formula \eqref{1.6a} for
$Z(\blam)$) as
           \begin{align}
            & C(I - Z(\blam)A^{*})^{-1}Z(\blam)Z(\bzeta)^{*} (I -
               A^{*}Z(\bzeta)^{*})^{-1}C^{*} + I_{\cY}  \notag \\
               & \qquad = C(I - Z(\blam) A)^{-1}(I - A^{*}
Z(\bzeta)^{*})^{-1}C^{*}+ S(\blam) S(\bzeta)^{*}.\label{eq}
\end{align}
It follows from the latter identity that the map
\begin{equation}
\label{defV}
V' \colon \begin{bmatrix} Z(\bzeta)^{*} (I - A^{*}
Z(\bzeta)^{*})^{-1} C^{*}  \\ I_{\cY}\end{bmatrix} y \mapsto
\begin{bmatrix} (I -A^{*}Z(\bzeta)^{*})^{-1}
C^{*} \\ S(\bzeta)^{*} \end{bmatrix} y
\end{equation}
can be extended by linearity and continuity to an isometry (still
denoted by $V'$) from the subspace
$${\mathcal D}_V : = \overline{\operatorname{span}} \left\{
\begin{bmatrix}Z(\bzeta)^{*}(I - A^{*}Z(\bzeta)^{*})^{-1}C^{*} \\
          I_{\cY} \end{bmatrix} y \colon \bzeta \in {\mathbb B}^{d} \text{
          and } y \in \cY \right\}
$$
onto the subspace ${\mathcal R}_V$ given
in \eqref{ranV}.  Note that the setting $\bzeta = 0 \in {\mathbb
B}^{d}$ in the formula $\sbm{Z(\bzeta)^{*}(I -
A^{*}Z(\bzeta)^{*})^{-1}C^{*}y \\ y }$ for a generic generator of
${\mathcal D}_{V}$ shows that $\{0\} \oplus \cY \subset {\mathcal
D}_{V}$ and hence we actually have
$ {\mathcal D}_{V} = {\mathcal D} \oplus \cY $
where ${\mathcal D}$ is defined as in \eqref{domV0}. Just as in the
proof of Proposition \ref{R:1.3},
setting $\bzeta = 0$ in the formula \eqref{defV} for the action
of $V'$ implies that
\begin{equation}
V' \colon \begin{bmatrix} 0 \\ y \end{bmatrix} \mapsto\begin{bmatrix}
C^{*} \\ S(0)^{*}\end{bmatrix} y\quad\mbox{for every} \; \; y\in\cY.
\label{vprime}
\end{equation}
Write $V'$ in the  block-matrix form $V' = \begin{bmatrix} A'_{V} &
B'_{V} \\ C'_{V} & D'_{V} \end{bmatrix}$ conformal with \eqref{defv}
and define $A_{V},B_{V},C_{V},D_{V}$ as in \eqref{abd} and \eqref{c}.
We conclude from \eqref{vprime} that $B'_{V}=C^*=B_V$, $D'_V=S(0)^*=D_V$.
Then \eqref{defV} implies that $C'_V$ satisfies
\begin{align*}
C_{V}' \left( Z(\bzeta)^{*}(I - A^{*}Z(\bzeta)^{*})^{-1}C^{*}y \right)
& = S(\bzeta)^{*}y - S(0)^{*}y \\
& = C_{V}\left( Z(\bzeta)^{*} (I -
A^{*}Z(\bzeta)^{*})^{-1} y \right)
\end{align*}
and hence $C_{V}' = C_{V}$.  Similarly,
\begin{align*}
A_{V}' \left(Z(\bzeta)^{*}(I - A^{*}Z(\bzeta)^{*})^{-1}C^{*}y \right)
       & = (I - A^{*}Z(\bzeta)^{*})^{-1}C^{*}y - C^{*}y \\&
       = A^{*} \left(Z(\bzeta)^{*}(I - A^{*}Z(\bzeta)^{*})^{-1}C^{*}y
        \right)
\end{align*}
and we conclude that $A_{V}' = A^{*}|_{\cD} = A_{V}$.
Thus, $V'=V$ and therefore $V$ is an isometry.

Conversely, if $V$ defined in \eqref{defv}--\eqref{c}
is isometric, then for two generic generators
$$
f=\begin{bmatrix} Z(\bzeta)^{*}(I
- A^{*}Z(\bzeta)^{*})^{-1} C^{*}y \\ y \end{bmatrix} \quad\text{and}\quad
g= \begin{bmatrix} Z(\blam)^{*}(I - A^{*}Z(\blam)^{*})^{-1}
           C^{*} y' \\ y' \end{bmatrix}
$$
in $\cD_{V}=\cD\oplus \cY$, we have
\begin{equation}
\left\langle f, \, g \right\rangle_{\cX^{d}
\oplus \cY}=\left\langle Vf, \, Vg \right\rangle_{\cX \oplus \cU}.
\label{eq1}
\end{equation}
Note that
           \begin{align*}
\left\langle f, \, g \right\rangle_{\cX^{d} \oplus \cY}&=
            \left\langle [ C(I-Z(\blam) A)^{-1}Z(\blam) Z(\bzeta)^{*}(I-A^{*}
            Z(\bzeta)^{*})^{-1}C^{*} + I ] y, \, y' \right
\rangle_{\cY} \notag \\
            & = \left \langle [ \langle \blam, \bzeta\rangle C(I -
            Z(\blam)A)^{-1} (I - A^{*}Z(\bzeta)^{*})^{-1}C^{*} + I ]y, y'
            \right\rangle_{\cY}
            \end{align*}
and
            \begin{align*}
\left\langle Vf, \, Vg \right\rangle_{\cX \oplus \cU}&=
\left\langle \begin{bmatrix} (I - A^{*}Z(\bzeta)^{*})^{-1}C^{*} \\
               S(\bzeta)^{*} \end{bmatrix} y, \,
               \begin{bmatrix} ( I - A^{*} Z(\blam)^{*})^{-1} C^{*} \\
               S(\blam)^{*} \end{bmatrix} y' \right \rangle_{\cX \oplus \cU}
               \notag \\
               & = \left\langle [ C (I - Z(\blam) A)^{-1}(I -
               A^{*}Z(\bzeta)^{*})^{-1}C^{*} + S(\blam) S(\bzeta)^{*}] y, y'
               \right  \rangle_{\cY}.
           \end{align*}
Substituting the two latter equalities into \eqref{eq1} and taking into
account that $y$ and $y'$ are arbitrary vectors in $\cY$, we get
\eqref{eq}, which is equivalent to \eqref{2.1a}.
\end{proof}

Proposition \ref{R:1.3} states that once a contractive realization
$\bU= \sbm{A & B \\ C & D}$ of $S$ is such that \eqref{2.1} holds, then
this realization is weakly coisometric. Our next result asserts that
equality \eqref{2.1} itself guarantees the existence of weakly
coisometric realizations for $S$ with preassigned $C$ and
$\bA=(A_1,\ldots,A_d)$.

\begin{theorem}  \label{T:CAtoS}
Suppose that a Schur-class function $S \in \cS_{d}(\cU, \cY)$ and
a contractive pair $(C, \bA)$ are such that \eqref{2.1} holds and let
$D:=S(0)$. Then there exist operators $B_{j} \colon \cU \to \cX$ for $j =
1, \dots, d$ so that the operator $\bU$ of the form \eqref{1.7a} is
weakly coisometric and $S$ can be realized as in \eqref{1.5a}.
              \end{theorem}
              \begin{proof} We are given $C$, $A$, $D = S(0)$
         and $S(\blam)$ for $\blam \in {\mathbb B}^{d}$ and seek $B \colon
         \cU \to \cX^d$ so that
         $$
          C (I - Z(\blam)A)^{-1} Z(\blam) B + D = S(\blam),
         $$
         or, in adjoint form with $\bzeta$ in place of $\blam$,
$$
B^{*} Z(\bzeta)^{*}(I-A^{*}Z(\bzeta)^{*})^{-1}C^{*}+D^{*}=S(\bzeta)^{*}.
$$
The latter equality is equivalent to
\begin{equation}
\label{want3}
\begin{bmatrix} A^{*} & C^{*} \\ B^{*} & D^{*} \end{bmatrix}
\begin{bmatrix} Z(\bzeta)^{*} (I - A^{*} Z(\bzeta)^{*})^{-1}
C^{*} \\ I \end{bmatrix} = \begin{bmatrix} (I - A^{*}
Z(\bzeta)^{*})^{-1} C^{*} \\ S(\bzeta)^{*} \end{bmatrix},
           \end{equation}
since the identity
$$
              A^{*}Z(\bzeta)^{*}(I - A^{*}Z(\bzeta)^{*})^{-1} C^{*} + C^{*} =
              (I - A^{*}Z(\bzeta)^{*})^{-1} C^{*}
$$
expressing equality of the top components in \eqref{want3} holds true
automatically. On the other hand, since by assumption \eqref{2.1}
holds, Lemma \ref{L:2.1} applies and the operator $V:
\cD\oplus\cY\to\cX\oplus\cU$
defined in \eqref{defv}--\eqref{c} is isometric and satisfies a similar
equality \eqref{defV} (with $V' = V$). It follows that
any choice of $B = \sbm{B_{1}
\\ \vdots \\ \\ B_{d}}$ such that
$\bU^*=\sbm{A^{*} & C^{*} \\ B^{*} & D^{*}}$ is a contractive
extension of $V$ from $\cD\oplus\cY$ to the whole of $\cX^d\oplus\cY$
gives rise to a weakly coisometric realization $\bU = \sbm{A & B \\ C
& D }$ for $S(\blam)$.

           Our completion problem (construction of $B$ subject to
           \eqref{want3} and that $\bU = \sbm{ A & B \\ C & D}$ be
           contractive) can now be reformulated as follows:  {\em Find an
           operator $B \colon \cU \to  \cX^{d}$ so that
           \begin{enumerate}
               \item the operator matrix
               $\begin{bmatrix} A^{*} & C^{*} \\ B^{*} & S(0)^{*}
               \end{bmatrix} \colon \begin{bmatrix} \cX^{d} \\ \cY
               \end{bmatrix} \to \begin{bmatrix} \cX \\ \cU \end{bmatrix}$
             is a contraction, and
             \item $B^{*}|_{\cD} = C_{V}$, where $C_V: \, \cD\to \cU$ is
given by \eqref{c}.
             \end{enumerate}}
\smallskip
             This is a contractive matrix-completion problem with linear
             side-constraint (2).  We  convert this problem to a standard
matrix-completion problem as follows. Let $ \cD^{\perp}:=\cX^d\ominus \cD$
and define operators
             $$
             T_{11} \colon \cD^{\perp} \to \cX,  \qquad
             T_{12} \colon   \cD \oplus \cY \to \cX, \qquad
             T_{22} \colon \cD \oplus \cY \to \cU
             $$
             by
             \begin{equation}  \label{defT}
             T_{11} = A^{*}|_{\cD^{\perp}}, \qquad
             T_{12} = \begin{bmatrix} A^{*}|_{\cD} & C^{*}
             \end{bmatrix}, \qquad
             T_{22} = \begin{bmatrix} C_{V} & S(0)^{*} \end{bmatrix}.
             \end{equation}
             Then our extension problem can be reformulated again as follows.
\begin{problem}
             Find an operator $X$ from $\cD^{\perp}$ to $\cU$ so that
             the block operator matrix
             \begin{equation}  \label{3.20}
	 \bU^{*} = \begin{bmatrix}  T_{11} & T_{12} \\ X & T_{22} \end{bmatrix}
             \colon \begin{bmatrix} \cD^{\perp} \\ \cD
             \oplus \cY \end{bmatrix} \to \begin{bmatrix} \cX \\ \cU
             \end{bmatrix}
             \end{equation}
             is a contraction.
\label{parrot}
\end{problem}
This is a standard matrix-completion problem handled by the
result of Parrott \cite{Parrott}:  {\em Problem \ref{parrot}
has a solution $X$ if and only if the obvious necessary conditions hold:
             \begin{equation} \label{necessary1}
          \left\| \begin{bmatrix} T_{11} & T_{12} \end{bmatrix}\right\| \le 1,
             \qquad
             \left\| \begin{bmatrix} T_{12} \\ T_{22} \end{bmatrix} \right\|
             \le 1
             \end{equation}}
Making use of the definitions of $T_{11}, T_{12},
             T_{22}$ from \eqref{defT}, we get more explicitly
             \begin{equation}  \label{necessary2}
	 \begin{bmatrix} T_{11} & T_{12} \end{bmatrix} =
             \begin{bmatrix} A^{*} & C^{*} \end{bmatrix}, \qquad
	 \begin{bmatrix} T_{12} \\ T_{22} \end{bmatrix} =
	     \begin{bmatrix} A^{*}|_{\cD} & C^{*} \\ C_{V} &
	     S(0)^{*} \end{bmatrix} = \begin{bmatrix} A_{V} & B_{V} \\
	     C_{V} & D_{V} \end{bmatrix} = V
              \end{equation}
              where we use the identification
              $$ \begin{bmatrix} \cD^{\perp} \\ \cD
              \oplus \cY \end{bmatrix} \cong  \begin{bmatrix} \cX^{d} \\
              \cY \end{bmatrix}$$
              in the first expression.
              Thus the first expression in \eqref{necessary2} is contractive
              by our assumption that $(C, \bA)$ is a contractive pair while
              the second expression collapses to $V$ which is isometric.  We
              conclude that the necessary conditions  \eqref{necessary1} are
              satisfied and hence, by the result of \cite{Parrott},
there exists a
              solution $X$ to  Problem \ref{parrot}. To complete the proof of
              Theorem \ref{T:CAtoS}, we set
              \begin{equation}  \label{3.14}
                B = \begin{bmatrix} X^{*} \\ C_{V}^{*} \end{bmatrix} \colon
                \cU \to \begin{bmatrix} \cD^{\perp} \\ \cD
                \end{bmatrix} \cong  \cX^{d}
              \end{equation}
              where $X$ is any solution of the matrix-completion problem
              \eqref{3.20}.  Note that the isometry property of $V$ then gives
              that the resulting colligation $\bU = \sbm{ A & B \\ C & D }$ is
              weakly coisometric.
              \end{proof}
\begin{remark}
{\rm   Every $X\in\cL(\cD^\perp,\cU)$ leading to a contractive
          (isometric or unitary) $\bU^*$ in \eqref{3.20}, gives rise via
formula \eqref{3.14} to a weakly coisometric (respectively, coisometric or
unitary) realization of $S$ of the form
\begin{equation}
\bU = \begin{bmatrix} A & ? \\ C & S(0)\end{bmatrix}.
\label{?}
\end{equation}}
\label{R:00}
\end{remark}
Applying well known descriptions \cite{Parrott, AG, dkw, T} of all $X$'s
solving contractive, isometric and unitary completion problems
\eqref{3.20} one can get all weakly coisometric, coisometric
or unitary realizations for $S$ of the form \eqref{?} as follows.
Let
          $T_{11},T_{12},T_{22}$ be as in \eqref{defT}.  Since
$\begin{bmatrix}
          T_{11} & T_{12} \end{bmatrix}$ is a contraction, there is a unique
          $G_{1} \colon \operatorname{\overline{Ran}} \,(I -
          T_{12}T_{12}^{*})^{1/2} \to \cD^{\perp}$ so that
          \begin{equation}  \label{3.25}
           G_{1}(I - T_{12} T_{12}^{*})^{1/2} = T_{11}^{*}, \qquad
            \operatorname{Ker}\, G_{1}^{*} = \operatorname{Ker}\, T_{11}.
          \end{equation}
          Since $\sbm{T_{12} \\ T_{22}} = V$ is an isometry, there exists a
          unique partial isometry
$$
G_{2} \colon \operatorname{\overline{Ran}}
          \, (I - T_{12}^{*}T_{12})^{1/2} = (\operatorname{Ker}\,
T_{22})^{\perp}
          \to \cU
$$
so that
          \begin{equation}  \label{3.26}
              G_{2}(I - T_{12}^{*} T_{12})^{1/2} = T_{22},\qquad
                 \operatorname{Ker}\, G_{2}^* = \operatorname{Ker}\,
T_{22}^*.
          \end{equation}
          The latter equality can be considered as the polar decomposition
of
          $T_{22}$. Note that
          \begin{equation}  \label{3.27}
          I_{\cU}-G_2G_2^*=P_{\operatorname{Ker}\, T_{22}^*}.
          \end{equation}
          (Here $P_{\operatorname{Ker}\, T_{22}^{*}}$ denotes the orthogonal
          projection onto $\operatorname{Ker}\, T_{22}^{*}$.)
          From the formula for $T_{22}$ in \eqref{defT} combined with the
          formula \eqref{c} for the action of $C_{V}$ on a generic
	generating vectors of $\cD$, we see that
          $$
          \operatorname{\overline{Ran}}\, T_{22} =
\operatorname{\overline{span}}
          \{ S(\bzeta)^{*}y \colon \; \bzeta \in {\mathbb B}^{d},\, y \in
               \cY \}
          $$
          and hence
\begin{equation}
\label{KerT22}
\operatorname{Ker}\, T_{22}^{*} = (\operatorname{\overline{Ran}} \,
T_{22})^{\perp} = \{ u \in \cU \colon S(\blam) u \equiv 0 \} =:\cU_{S}^{0}.
\end{equation}
          \begin{theorem}  \label{T:3.5}
          Suppose that $S \in \cS_{d}(\cU, \cY)$ and a
          contractive pair $(C,\bA)$ are such that $K_{S}(\blam, \bzeta) =
          K_{C,\bA}(\blam, \bzeta)$.  Let $\cD \subset \cX^{d}$, $C_V$,
          $T_{11}$, $T_{12}$ and $T_{22}$ be as in \eqref{domV0},
          \eqref{c}, \eqref{defT} with $G_{1}$, $G_{2}$
          constructed as in \eqref{3.25}, \eqref{3.26} and the subspace
          $\cU^{0}_S$ given as in \eqref{KerT22}.  Then:
          \begin{enumerate}
          \item A realization ${\bf U}=\sbm {A & B \\ C & S(0)}$ of $S$
          is weakly coisometric if and only if $B$ is of the form
\begin{equation}  \label{desB}
                B = \begin{bmatrix} X^{*} \\ C_{V}^{*} \end{bmatrix}\quad
\mbox{where}\quad
X=-G_2T_{12}^*G_1^*+Q(I_{\cD^\perp}-
                     G_1G_1^*)^{\frac{1}{2}}
              \end{equation}
and where $Q\colon\operatorname{\overline{Ran}}\,
(I_{\cD^\perp}-G_1G_1^*)^{1/2}\to \cU_S^{0}$ is a contraction.
          \item $S$ admits a coisometric realization $\bU$ of the form
               \eqref{?} if and only if
          \begin{equation}
          \dim \operatorname{\overline{Ran}}\,
          (I_{\cD^\perp}-G_1G_1^*)^{1/2}\le \dim \, \cU_S^{0}.
          \label{3.29}
          \end{equation}
In this case, a realization ${\bf U}=\sbm {A & B \\ C & S(0)}$ of $S$
is coisometric if and only if $B$ is of the form \eqref{desB} for
some isometric $Q\colon\operatorname{\overline{Ran}}\,
(I_{\cD^\perp}-G_1G_1^*)^{1/2}\to \cU_S^{0}$.
          \item $S$ admits a unitary realization $\bU$ of the form \eqref{?}
          if and only if $(C, \bA)$ is an isometric pair, i.e.
         \begin{equation}
              A_{1}^{*}A_{1} + \cdots + A_{d}^{*} A_{d} + C^{*}C =
              I_{\cX},
          \label{3.29a}
          \end{equation}
           and
          \begin{equation}
          \dim \, \left(\operatorname{Ker}\, A^*\cap
\cD^\perp\right)= \dim \,
          \cU_S^{0}.
          \label{3.30}
          \end{equation}
In this case, a realization ${\bf U}=\sbm {A & B \\ C & S(0)}$ of $S$
is unitary if and only if $B$ is of the form \eqref{desB} for
some unitary $Q\colon\operatorname{\overline{Ran}}\,
(I_{\cD^\perp}-G_1G_1^*)^{1/2}\to \cU_S^{0}$.
          \end{enumerate}
          \end{theorem}
          \begin{proof}
             Problem \ref{parrot} is  equivalent to the following positive
completion problem: {\em find $X$  such that}
                \begin{equation}
                \begin{bmatrix}I & 0 & T_{11}^* & X^* \\
                0 & I & T_{12}^* & T_{22}^*\\ T_{11} & T_{12} & I & 0 \\
                X & T_{22} & 0 & I\end{bmatrix}\ge 0.
                \label{3.31}
                \end{equation}
                Substituting expressions \eqref{3.25} and \eqref{3.26}
                for $T_{11}^{*}$ and $T_{22}$ into
                \eqref{3.31} and taking the Schur complement to the
                principal
                (positive semidefinite) block $\sbm{I & T_{12}^* \\ T_{12} &
                I}$ we get (upon invoking \eqref{3.27} and \eqref{KerT22})
                \begin{equation}
                \begin{bmatrix}
                I_{\cD^\perp}-G_1G_1^* & X^*+G_1T_{12}G_2^* \\
                X+G_2T_{12}^*G_1^* & P_{\cU_S^{0}}\end{bmatrix}\ge 0
                \label{3.32}
                \end{equation}
                which is equivalent to \eqref{3.31}. It follows from
                \eqref{3.32} and \eqref{3.27} that $X$ is a solution of the
                contractive completion problem \eqref{3.20}
                (and therefore it leads via formula \eqref{3.14} to a
                weakly coisometric realization of $S$) if and only if it is
                of the form
                \begin{equation}
                X=-G_2T_{12}^*G_1^*+Q(I_{\cD^\perp}-
                     G_1G_1^*)^{\frac{1}{2}}
                \label{3.33}
                \end{equation}
                for some contraction $Q \colon\operatorname{\overline{Ran}}\,
                (I_{\cD^\perp}-G_1G_1^*)^{1/2}\to \cU_S^{0}$ which
                on account of Remark \ref{R:00} completes
                the proof of the first statement in the theorem.

Note that a contractive $\bU^*$ of the form \eqref{3.20} is an
isometry if  and only if
\begin{equation}
T_{11}^*T_{11}+X^*X=I_{\cD^\perp}.
\label{3.33a}
\end{equation}
To simplify the latter relation we need the following two equalities:
\begin{equation}
G_2^*Q(I_{\cD^\perp}-G_1G_1^*)^{\frac{1}{2}}=0\quad\mbox{and}\quad
G_1T_{12}\vert_{\operatorname{Ker}\, G_2}=0.
\label{3.33b}
\end{equation}
The first equality holds true since $ \operatorname{Ran}\, Q\subset
\cU_S^{0}=\operatorname{Ker}\, T_{22}^*=\operatorname{Ker}\, G_2^*$, by
\eqref{KerT22} and \eqref{3.26}. To verify the second equality,
take a vector $x\in\operatorname{Ker}\, G_2$ in the form
$$
x=(I - T_{12}^{*} T_{12})^{1/2}y\quad\mbox{where}\quad
y\in\operatorname{Ker} \, T_{22}.
$$
Then, by  \eqref{3.25},
\begin{equation}
G_1T_{12}x=G_1T_{12}(I - T_{12}^{*} T_{12})^{1/2}y=
G_1(I - T_{12}T_{12}^*)^{1/2}T_{12}y=T_{11}^*T_{12}y.
\label{3.33c}
\end{equation}
Since $\begin{bmatrix} T_{11} & T_{12} \end{bmatrix}$ is a contraction,
$$
\|T_{11}^*T_{12}y\|^2+\|T_{12}^*T_{12}y\|^2\le \|T_{12}y\|^2
$$
and since  $\begin{bmatrix} T_{12} \\ T_{22} \end{bmatrix}$ is an
isometry and $T_{22}y=0$, we have
$$
        \|T_{12}y\|=\|T_{12}^*T_{12}y\|=\|y\|.
$$
Combining the two latter relations we conclude that $T_{11}^*T_{12}y=0$
and now the second relation in \eqref{3.33b} follows from
\eqref{3.33c}. Making use of \eqref{3.33} and of the first relation in
\eqref{3.33b}, we get
$$
X^*X=(I-G_1G_1^*)^{\frac{1}{2}}Q^*Q(I-G_1G_1^*)^{\frac{1}{2}}
+G_1T_{12}G_2^*G_2T_{12}^*G_1^*
$$
which being substituted along with \eqref{3.25} into \eqref{3.33a}
allows us to write \eqref{3.33a} equivalently as
$$
(I-G_1G_1^*)^{\frac{1}{2}}(I-Q^*Q)(I-G_1G_1^*)^{\frac{1}{2}}
=-G_1T_{12}(I-G_2^*G_2)T_{12}^*G_1^*.
$$
Since $I-G_2^*G_2$ is equal to the orthogonal projection onto
$\operatorname{Ker}\, G_2$, the expression on the right hand side equals
zero and thus, \eqref{3.33a} is equivalent to
$$
(I_{\cD^\perp}-G_1G_1^*)^{\frac{1}{2}}
(I-Q^*Q)(I_{\cD^\perp}-G_1G_1^*)^{\frac{1}{2}}=0
$$
which means that $Q$ is isometric. The latter may occur if and only if
condition \eqref{3.29} holds. This completes the proof of the second
statement in the theorem.

Finally, for $\bU^*$ to be unitary it is
necessary that  $\begin{bmatrix} T_{11} & T_{12}\end{bmatrix}$ is a
coisometry, which
on account of \eqref{necessary2} can be written as
$A^*A+C^*C=I_{\cX}$ and is equivalent to \eqref{3.29a}. In this case
the operator $G_1$  defined in \eqref{3.25} is a partial isometry and
                $$
                I_{\cD^\perp}-G_1G_1^*=P_{{\rm Ker} \, T_{11}}.
                $$
Then the parametrization formula \eqref{3.33} for all solutions $X$ of
the  contractive completion problem takes the form
                \begin{equation}
                X=G_2T_{12}^*G_1^*+Q
                \label{3.34}
                \end{equation}
                where $Q: \, {\rm Ker} \, T_{11}\to \cU_S^{0}$ is a
contraction.
                A contraction $\bU^*$ of the form \eqref{3.20} is an
                     isometry if and only if
                $$
                T_{11}^*T_{11}+X^*X=I_{\cD^\perp}\quad\mbox{and}
\quad XX^*+T_{22}T_{22}^*=I_{\cU}.
                $$
Substituting \eqref{3.25} and  \eqref{3.34} into the latter
equalities we write them equivalently as
$$
I_{{\rm Ker} \, T_{11}}-Q^*Q=0\quad\mbox{and}\quad
I_{\cU_S^0}-QQ^*=0
$$
which means that $Q$ must be unitary. The latter may occur if and
only if
                $$
                \dim \, {\rm Ker} \, T_{11}= \dim \, \cU_S^0.
                $$
                Since $\operatorname{Ker} \, T_{11}=\operatorname{Ker}\,
                A^*|_{\cD^\perp}=
                \operatorname{Ker} \, A^*\cap \cD^\perp$, the
                     last condition is equivalent to \eqref{3.30}.
\end{proof}
As a corollary we obtain the following uniqueness result.
\begin{corollary}
\label{C:3.5}
Suppose that $S \in \cS_{d}(\cU, \cY)$ and a
contractive pair $(C,\bA)$ are such that $K_{S}(\blam, \bzeta) =
K_{C,\bA}(\blam, \bzeta)$.  Let $\cD \subset \cX^{d}$,
$T_{11}$, $T_{12}$ be as in \eqref{domV0}, \eqref{defT} with $G_{1}$,
	constructed as in \eqref{3.25}, and the subspace
	$\cU^{0}_S$ given as in \eqref{KerT22}.  Then:
	\begin{enumerate}
	\item $S$ admits a unique weakly coisometric realization $\bU$ of the
	form \eqref{?} if and only if either
	\begin{equation}
	G_1G_1^*=I_{\cD^\perp}\quad \text{or}\quad
               \cU_S^{0}=\{0\}.
	\label{3.28}
	\end{equation}
          \item If $G_1G_1^*=I_{\cD^\perp}$, then this unique realization
          is also coisometric and it is unitary if  $(C, \bA)$ is an
          isometric pair and
          both conditions in \eqref{3.28} are satisfied.
\item In either case, this unique realization is obtained via formula
	\eqref{desB} applied to  $X=-G_2T_{12}^*G_1^*$.
	\end{enumerate}
	\end{corollary}
The second condition in \eqref{3.28} is much easier to be verified.
We display uniqueness caused by this condition as a separate statement.
	  \begin{corollary} \label{C:3.4}  Let $S \in \cS_d(\cU, \cY)$
                 and let $(C, \bA)$ be a contractive pair
	  such that $K_{S}(\blam, \bzeta) = K_{C, \bA}(\blam, \bzeta)$.
	  Suppose that $\cU^{0}_{S} = \{0\}$, i.e., that
	  \begin{equation}  \label{3.35}
	      S(\blam) u \equiv 0 \Longrightarrow u = 0.
	   \end{equation}
	   Then $S$ admits a unique weakly coisometric realization $\bU$ of
	   the form \eqref{1.5a} consistent with the preassigned
	   choice of output pair $(C, \bA)$.  Moreover:
	   \begin{enumerate}
	   \item This realization is coisometric if and only if
	   $G_1G_1^*=I_{\cD^\perp}$, where $G_1$ is defined in \eqref{3.25}.
	   \item This realization is unitary if and only if
                 $(C, \bA)$ is an isometric pair
                and $\operatorname{Ker}\, A^*\cap \cD^\perp=\{0\}$.
	   \end{enumerate}
	   \end{corollary}

	   The case when $S$ satisfies condition \eqref{3.35} is generic in the
	   following sense: if the subspace $\cU_S^0$ is not trivial, we
	   represent $\cU$ as $(\cU_S^0)^\perp\oplus\cU_S^0$ and
	   write $S(\blam)$ with respect to this decomposition as
	   $$
	   S(\blam)=\begin{bmatrix}\widetilde{S}(\blam) & 0 \end{bmatrix}.
	   $$
	   Then $\widetilde{S}\in\cS_d((\cU_S^0)^\perp,\cY)$ satisfies the
	   condition
	   \eqref{3.35} and besides, $K_{C, \bA}(\blam, \bzeta) =
	   K_{S}(\blam, \bzeta) = K_{\widetilde S}(\blam, \bzeta)$.
	   Suppose that we are given a contractive pair $(C, \bA)$ such
	   that $K_{S}(\blam, \bzeta) = K_{C, \bA}(\blam, \bzeta)$ and we let
$$
\widetilde{S}(\blam)=\widetilde{D}+C(I-Z(\blam)A)^{-1}Z(\blam)\widetilde{B}
$$
	   be the unique weakly coisometric realization of $\widetilde{S}$
	   consistent with $(C, \bA)$ and  $\widetilde D=\widetilde S(0)$.
	   Then every weakly coisometric realization for $S$
	   consistent with $(C, \bA)$ and  $\widetilde S(0)$ is of the form
	   \eqref{1.5a} with
	   $$
	   D=\begin{bmatrix}\widetilde{D} & 0 \end{bmatrix}\quad\mbox{and}\quad
	   B=\begin{bmatrix}\widetilde{B}& B^0\end{bmatrix}
	   $$
	   where $B^0: \, \cU_S^0\to \cX^d$ is an operator
                  subject to the sole constraint that the operator
	   \begin{equation}
	   \bU=\begin{bmatrix} A & \widetilde{B}& B^0 \\
	   C & \widetilde{D} & 0 \end{bmatrix}
	   \label{3.36}
	   \end{equation}
	   be a contraction. This operator $B^0$ is responsible
                  for nonuniqueness of weakly coisometric realizations
                  compatible with a given contractive pair $(C, \bA)$; it is
                  also clear that if $\dim \, \cU_S^0$ is large enough,
                  $\bU$ of the form  \eqref{3.36} can be arranged to be
                  coisometric. We can look at this from another point of
	   view as follows.

	   \begin{proposition}
	   If $S\in\cS_d(\cU,\cY)$ admits a weakly coisometric
	   realization, then there exists a Hilbert space ${\mathcal F}$ and
	   a partial isometry $W\colon \, \cF\to \cU$ so that the function
	   $S_W(z)=S(z)W\in\cS_d(\cF,\cY)$ admits  a coisometric
	   realization. If in addition condition \eqref{3.29a} is
satisfied, then $\cF$ and $W$ can be chosen so that $S_W$ admits a unitary
	   realization.
	   \label{R:3.10}
	   \end{proposition}
	   \begin{proof}
	   It suffices to pick $\cF=(\cU_S^0)^\perp\oplus
	   \operatorname{\overline{Ran}}\,
	   (I_{\cD^\perp}-G_1G_1^*)^{1/2}$ and to define
the partial isometry
	   $W \colon  \cF \to \cU$ by $Wf=f$ if $f\in(\cU_S^0)^\perp$ and
	   $Wf=0$ if $f\in \operatorname{\overline{Ran}}\,
(I_{\cD^\perp}-G_1G_1^*)^{1/2}$. \end{proof}

The analysis in the proofs of Theorems \ref{T:CAtoS} and \ref{T:3.5}
can be slightly modified to get a description of all Schur class
representers of a contractive pair $(C,\bA)$.
\begin{theorem}  \label{T:repr}
Let $(C, \bA)$ be a contractive pair with $C\in\cL(\cX,\cY)$, let $\cD$ be
the subspace of $\cX^d$  given by \eqref{domV0} and let
\begin{equation}
T:=\begin{bmatrix} A^*\vert_{\cD} & C^*\end{bmatrix}: \; \cD\oplus
\cY\to\cX.
\label{new1}
\end{equation}
\begin{enumerate}
\item Given a Hilbert space $\cU$, there exists an $S\in\cS_d(\cU,\cY)$
such that
\begin{equation}
K_{S}(\blam, \bzeta) = K_{C, \bA}(\blam, \bzeta)
\label{new2}
\end{equation}
if and only if
\begin{equation}
\dim \, \cU\ge \dim \,
\overline{\operatorname{Ran}}\, (I-T^*T)^{\frac{1}{2}}.
\label{new3}
\end{equation}
\item If \eqref{new3} is satisfied, then all $S\in\cS_d(\cU,\cY)$
for which \eqref{new2} holds are described by the formula
\begin{equation}
S(\blam)=\begin{bmatrix} C(I-Z(\blam)A)^{-1} &
I_\cY\end{bmatrix}(I-T^*T)^{\frac{1}{2}}G^*
\label{new3a}
\end{equation}
where $G$ is an isometry from
$\overline{\operatorname{Ran}}\, (I-T^*T)^{\frac{1}{2}}$ onto
$\operatorname{Ran}\, G\subset\cU$.
\item If $\dim \, \cU=\dim \,
\overline{\operatorname{Ran}}\,(I-T^*T)^{\frac{1}{2}}$, then the function
$S\in\cS_d(\cU,\cY)$ such that \eqref{new2} holds is defined uniquely up
to a constant unitary factor on the right.
\end{enumerate}
\end{theorem}
\begin{proof} By Lemma \ref{L:2.1}, if there is an $S\in\cS_d(\cU,\cY)$
such that \eqref{new2} holds, then the operator $V$ defined
\eqref{defv}--\eqref{c} is an isometry. It is readily seen from
\eqref{abd}, \eqref{c} and \eqref{new1} that the top block row in $V$ is
equal to $T$ while the bottom block row
\begin{equation}
\widetilde{T}:=\begin{bmatrix} C_V & D_V\end{bmatrix}: \;
\cD\oplus \cY\to\cU
\label{new4}
\end{equation}
depends on $S(\blam)$ and is not specified in the conditions of the
theorem. Thus, a necessary condition for an $S\in\cS_d(\cU,\cY)$ to
exist so that \eqref{new2} holds is that there exists $\widetilde{T}: \,
\cD\oplus \cY\to\cU$ such that  the operator
\begin{equation}
V=\begin{bmatrix} T \\ \widetilde{T}\end{bmatrix}: \;
\begin{bmatrix} \cD\\  \cY\end{bmatrix}\to \begin{bmatrix} \cX\\
\cU\end{bmatrix}
\label{new5}
\end{equation}
is isometric. The latter is true if and only if the condition
\eqref{new3} is satisfied (which proves the necessity part in Statement
(1) of the theorem) and every such $\widetilde{T}$ is necessarily of
the form
\begin{equation}
\widetilde{T}=G(I-T^*T)^{\frac{1}{2}}
\label{new6}
\end{equation}
where $G$ is an isometry from
$\overline{\operatorname{Ran}}\,(I-T^*T)^{\frac{1}{2}}$ onto
$\operatorname{Ran}\, G\subset\cU$. The equality
\begin{equation}
S(\bzeta)^*y=\widetilde{T}\begin{bmatrix}
Z(\bzeta)^*(I_{\cX}-A^{*}Z(\bzeta)^{*})^{-1} C^{*}\\ I_{\cY}\end{bmatrix}y
\qquad (\bzeta\in\B^d; \; y\in\cY)
\label{new7}
\end{equation}
defines an $\cL(\cU,\cY)$-valued function $S(\bzeta)$ pointwise. By
setting $\bzeta=0$ in \eqref{new7} we get
$$
S(0)^*y=\widetilde{T}\begin{bmatrix} 0\\ I_{\cY}\end{bmatrix}y
\qquad (\bzeta\in\B^d; \; y\in\cY),
$$
and therefore, the block entry $D_V$ in \eqref{new4} is equal to
$S(0)^*$. Then it follows from \eqref{new7} that the  block entry
$C_V$ in \eqref{new4} is defined explicitly as in the formula
\eqref{c}. Thus, the isometry $V$ in \eqref{new5} coincides with that
in  \eqref{defv}--\eqref{c}. Then we apply Lemma \ref{L:2.1} to
conclude that \eqref{new2} holds for $S$ defined in \eqref{new7}
and in particular, that this $S$ belongs to
$\cS_{d}(\cU, \cY)$. This completes the proof of Statement (1).

Since every representer $S$ gives rise to an isometric extension $V$
of $T$ as in \eqref{new5} and since \eqref{new6} is the general
formula for the bottom component of $V$, it follows that the formula
\eqref{new7} gives a parametrization of all representers
$S\in\cS_d(\cU,\cY)$. Replacing $\widetilde{T}$ in \eqref{new7} by
its expression \eqref{new6} and taking into account that $y\in\cY$ is
arbitrary, we get
$$
S(\bzeta)^*y=G(I-T^*T)^{\frac{1}{2}}\begin{bmatrix}
Z(\bzeta)^*(I_{\cX}-A^{*}Z(\bzeta)^{*})^{-1} C^{*}\\ I\end{bmatrix}.
$$
Taking adjoints we arrive at \eqref{new3a}. The last statement of the
theorem
now is self-evident, since under the assumption that $\dim \, \cU=\dim  \,
\overline{\operatorname{Ran}}\,(I-T^*T)^{\frac{1}{2}}$, the operator $G$
is unitary.
\end{proof}

\section{Generalized functional-model realizations}
\label{S:com}

Constructing a weakly coisometric realization for a given $S \in
\cS_{d}(\cU, \cY)$ is not an issue: by Theorem \ref{T:BTV}, every $S \in
\cS_{d}(\cU, \cY)$ admits even a unitary realization. However, the pair
$(C, \bA)$  for a weakly coisometric realization can be constructed in a
certain canonical way.

         \begin{theorem} \label{T:noncomreal}
Let $S \in \cS_{d}(\cU, \cY)$ and let $\cH(K_{S})$ be the associated de
Branges-Rovnyak space. Then:
\begin{enumerate}
\item There exist bounded operators $A_{j} \colon \cH(K_{S})
\to \cH(K_{S})$ such that
\begin{equation}
f(\blam) - f(0)  = \sum_{j=1}^{d} \lam_{j} (A_{j}f)(\blam)
\quad\mbox{for every $f \in \cH(K_{S})$ and $\blam \in {\mathbb B}^{d}$},
\label{dop1}
\end{equation}
and
\begin{equation}
        \sum_{j=1}^{d} \|A_{j}f\|^{2}_{\cH(K_{S})}  \le
             \|f\|^{2}_{\cH(K_{S})} - \| f(0)\|^{2}_{\cY}.
\label{dop2}
\end{equation}
\item  There is a weakly coisometric realization \eqref{1.5a} for $S$
with state space $\cX$ equal to $\cH(K_{S})$ with the state operators
$A_1,\ldots, A_d$ from part (1) and the operator $C \colon  \cH(K_S)\to \cY$
defined by
\begin{equation}
Cf=f(0)  \quad\text{for all}\quad f\in\cH(K_S).
\label{dop8}
\end{equation}
\end{enumerate}
\end{theorem}
\begin{proof} Since every $S \in \cS_{d}(\cU, \cY)$ admits a weakly
coisometric realization, the associated space $\cH(K_S)$ can be identified
as the range space of the observability operator for some contractive
pair. Then part (1) of the theorem follows from Theorem \ref{T:3-1.2nc}.
Now let us assume that relations \eqref{dop1} and
\eqref{dop2} hold and that $C$ is defined as in \eqref{dop8}. Then
\eqref{dop2} says that the pair $(C, \bA)$ is contractive.
Iteration of \eqref{dop1} says that, for each $f \in\cH(K_S)$,
\begin{align*}
f(\blam) = & \sum_{j_{1}=1}^{d} \lambda_{j_{1}}\left[
(A_{j_{1}}f)(0)+ \sum_{j_{2}=1}^{d} \lambda_{j_{2}} \left[(A_{j_{2}}
A_{j_{1}}f)(0) + \sum_{j_{3}=1}^{d}\lambda_{j_{3}}\left[(A_{j_{3}}
               A_{j_{2}}  A_{j_{1}}f(0) +  \right.\right.\right.\\
               & \qquad \left.\left.\cdots + \sum_{j_{k}=1}^{d}
              \lambda_{j_{k}}\left[(A_{j_{k}}
              \cdots  A_{j_{2}}A_{j_{1}}f)(0) + \cdots \right] \cdots
              ]\right]\right].
               \end{align*}
               This unravels to the tautology
\begin{equation}  \label{reproduce}
           f(\blam) = C (I - Z(\blam) A)^{-1}f \quad\text{for
           all}\quad f\in\cH(K_{S}).
         \end{equation}
         Hence, by the reproducing property of $K_{S}$,
         for any $\bzeta \in {\mathbb B}^{d}$, $y \in \cY$ and $f \in
         \cH(K_{S})$, we have
         \begin{align*}
             \langle f, \; K_{S}(\cdot, \bzeta) y \rangle_{\cH(K_{S})} & =
             \langle f(\bzeta), \; y \rangle_{\cY} \\
             & = \langle C(I - Z(\bzeta) A)^{-1}f, y \rangle_{\cY} \\
             & = \langle f, \; (I - A^{*} Z(\bzeta)^{*})^{-1}C^{*} y
             \rangle_{\cH(K_{S})}
         \end{align*}
         and we conclude that
         \begin{equation}  \label{H-canonical}
         K_{S}(\cdot, \bzeta) y = (I - A^{*} Z(\bzeta)^{*})^{-1} C^{*} y.
\end{equation}
Hence, for all $\blam, \bzeta \in
         {\mathbb B}^{d}$ and $y,y' \in \cY$ we have
         \begin{align}
             \langle K_{S}(\blam, \bzeta) y, \; y' \rangle_{\cY} & =
             \langle K_{S}(\cdot, \bzeta) y, \; K_{S}( \cdot, \blam) y'
             \rangle_{\cH(K_{S})} \notag\\
             & = \langle (I - A^{*}Z(\bzeta)^{*})^{-1}C^{*}y, \; (I -
             A^{*}Z(\blam)^{*})^{-1}C^{*}y' \rangle_{\cH(K_{S})}\notag \\
             & = \langle C(I - Z(\blam) A)^{-1}(I -
             A^{*}Z(\bzeta)^{*})^{-1}C^{*}y,\;  y' \rangle_{\cY}\notag \\
& = \langle K_{C,\bA}(\blam, \bzeta) y, \; y' \rangle_{\cY}\label{aga}
         \end{align}
         from which we conclude that
$K_{S}(\blam, \bzeta) = K_{C, \bA}(\blam, \bzeta)$.
         It now follows from Theorem \ref{T:CAtoS} that there is a
         choice of $B_{j} \colon \cU \to \cH(K_S)$ with $\bU = \sbm{
         A & B \\ C & D } \colon \cH(K_S) \oplus \cU \to \cH(K_S)^{d} 
\oplus \cY$
         weakly coisometric so that $S(\blam) = D + C(I - Z(\blam)A)^{-1}
         Z(\blam)B$.  This completes the proof.\end{proof}
Equality \eqref{dop1} means that the operator tuple
${\bf A}=(A_1,\ldots,A_d)$ solves {\em the Gleason problem} \cite{gleason}
for $\cH(K_{S})$. Let us say that ${\bf A}$ {\em is a
contractive solution of the Gleason problem} if in addition relation
\eqref{dop2} holds for every $f\in\cH(K_S)$ or, equivalently, if the pair
$( C,{\bf A})$ is contractive where $C: \, \cH(K_S)\to\cY$ is defined as
in \eqref{dop8}. Theorem \ref{T:noncomreal} shows that any contractive
solution ${\mathbf A} = (A_{1},\dots, A_{d})$ of the Gleason problem for
$\cH(K_{S})$ gives rise to a weakly coisometric realization for
$S \in\cS_{d}(\cU, \cY)$ (not unique, in general). Let us call any such
weakly coisometric realization  a {\em generalized functional-model}
realization of $S(\blam)$.
A consequence of  formula \eqref{reproduce} is that {\em any generalized
functional-model realization of $S$ is observable.}

Note also that any contractive realization
\begin{equation}
{\mathbf U} = \begin{bmatrix} A & B \\ C & D \end{bmatrix} \colon
\begin{bmatrix} \cH(K_{S}) \\ \cU \end{bmatrix} \to \begin{bmatrix}
          \cH(K_{S})^{d} \\ \cY \end{bmatrix}
\label{note}
\end{equation}
with $C$ given as in \eqref{dop8} and the state space tuple
$(A_1,\ldots,A_d)$ a contractive solution to the Gleason problem on
$\cH(S)$ is automatically weakly coisometric (i.e., a generalized
functional-model realization), as follows from
calculation \eqref{aga} and Proposition \ref{R:1.3}.

For a generalized functional-model realization, we have the
following explicit formulas for the characters appearing in Lemma
\ref{L:2.1}.
\begin{proposition}  \label{P:canonical-guys}
Suppose that ${\mathbf U}$ of the form \eqref{note}
is a generalized functional model realization for the Schur-class
function $S \in {\mathcal S}_{d}(\cU, \cY)$ and that the spaces
${\mathcal D}$ and ${\mathcal R}_{V}$ are defined as in \eqref{domV0}
and \eqref{ranV}.  Then the spaces ${\mathcal D}$, ${\mathcal
D}^{\perp} = \cH(K_{S})^{d} \ominus {\mathcal D}$, ${\mathcal R}_{V}$
and ${\mathcal R}_{V}^{\perp} = \cH(K_{S}) \ominus {\mathcal R}_{V}$
can be described in the following explicit functional forms:
          \begin{align}
              {\mathcal D} & = \operatorname{\overline{span}} \{
              Z(\bzeta)^{*} K_{S}(\cdot, \bzeta) y \colon \bzeta \in
{\mathbb
              B}^{d},\, y \in \cY \}, \notag \\
              \cR_{V} & = \operatorname{\overline{span}}
              \left\{ \begin{bmatrix} K_{S}(\cdot, \bzeta) y, \\
             S(\bzeta)^{*} y \end{bmatrix} \colon \bzeta \in {\mathbb
              B}^{d}, \, y \in \cY \right\}, \notag \\
               {\mathcal D}^{\perp} & = \{ h  \in \cH(K_{S})^{d} \colon
             Z(\blam) h(\blam) \equiv 0 \}, \label{canonical3} \\
             \cR_{V}^{\perp} & = \left\{ \begin{bmatrix} h \\ u
             \end{bmatrix} \in \begin{bmatrix} \cH(K_{S}) \\ \cU
             \end{bmatrix} \colon h(\blam) + S(\blam) u \equiv 0 \right\}.
             \label{canonical2}
          \end{align}
         \end{proposition}
          \begin{proof} Substituting \eqref{H-canonical} into  \eqref{domV0}
and
\eqref{ranV} gives the two first of the four
representations above. Given the formula for $\cD$, the formula
          for $\cD^{\perp}$ follows via a standard calculation in
          reproducing kernel Hilbert spaces using the reproducing
          kernel property: indeed $f \in \cH(K_{S})^{d} \ominus \cD$
          means that
          $$
          0 = \langle f, Z(\bzeta)^{*} K_{S}(\cdot, \bzeta) y
          \rangle_{\cH(K_{S})^{d}} =
          \langle Z(\bzeta) f, K_{S}(\cdot, \bzeta) \rangle_{\cH(K_{S})} =
          \langle Z(\bzeta) f(\bzeta), y \rangle_{\cY}
          $$
          holds for every $\bzeta \in {\mathbb B}^{d}$ and $y \in \cY$
          forcing $Z(\blam) f(\blam) \equiv  0$.
          \end{proof}
\begin{remark}
\label{R:collapse}
{\rm Note that in case $d=1$, Theorem \ref{T:noncomreal}
collapses to Theorem \ref{T:clas-dBRreal}. Indeed, in this case
equality \eqref{dop1} reads
$$
f(\lambda)-f(0)=\lambda (Af)(\lambda)
$$
and completely defines the operator $A$ as in formula
\eqref{clas-dBRrealac}. By \eqref{canonical3}, ${\mathcal D}^\perp=\{0\}$
and hence ${\mathcal D}=\cH(K_S)$. Therefore any weakly
coisometric realization is automatically coisometric.
On account of \eqref{H-canonical},
formula \eqref{c} for $C_V:{\mathcal D}=\cH(K_S)\to\cU$ takes the form
$$
C_{V} \colon \; \bar{\zeta}K_{S}(\cdot, \bzeta) y \
\mapsto (S(\zeta)^{*} - S(0)^{*}) y\quad \mbox{for}\; \;
\zeta\in{\mathbb D}\; \; y\in\cY
$$
and the formula for its adjoint $C_{V}^*: \, \cU\to \cH(K_S)$,
\begin{equation}
\label{cv*}
C_{V}^* \colon \; u\to \frac{S(\zeta) - S(0)}{\zeta} \, u,
\end{equation}
follows from equalities
\begin{eqnarray*}
\langle  (C_{V}^*u)(\zeta), \, y\rangle_{\cY}&=&
\langle  C_{V}^*u, \, K_{S}(\cdot, \zeta) y\rangle_{\cH(K_S)}\\
&=&
\langle u, \, C_{V} K_{S}(\cdot, \zeta)y\rangle_{\cU}\\
&=&\left\langle u, \, \frac{S(\zeta)^{*} - S(0)^{*}}{\bar{\zeta}}
y\right\rangle_{\cU}
=\left\langle \frac{S(\zeta) - S(0)}{\zeta} \, u, y\right\rangle_{\cY}.
\end{eqnarray*}
Since ${\mathcal D}^\perp = \{0\}$, formula \eqref{3.14} gives that the only
$B$ such that $\bU=\sbm{A & B \\ C & D}$ is coisometric is $B=C_{V}^*$.
By \eqref{cv*}, this $B$ is the same as in \eqref{clas-dBRrealac}.}
\end{remark}
We next present the result concerning the
universality of generalized functional-model realizations among
weakly coisometric realizations. We say that two colligations
$$
{\mathbf U} = \begin{bmatrix}A & B \\ C & D\end{bmatrix}
\colon\cX \oplus \cU \to \cX^{d} \oplus \cY\quad\mbox{and}\quad
\widetilde{{\mathbf U}}
= \begin{bmatrix}\widetilde{A}& \widetilde{B} \\ \widetilde{C} &
\widetilde{D}\end{bmatrix}
\colon\widetilde{\cX} \oplus \cU \to \widetilde{\cX}^{d}
\oplus \cY
$$
are {\em unitarily equivalent} if
there is a unitary operator $U \colon \cX \to \widetilde{\cX}$ such that
$$
\begin{bmatrix} \oplus_{k=1}^{d} U & 0 \\ 0 & I_{\cY} \end{bmatrix}
\begin{bmatrix} A & B \\ C & D \end{bmatrix} =
\begin{bmatrix}\widetilde{A}&
\widetilde{B} \\ \widetilde{C} &
\widetilde{D} \end{bmatrix} \begin{bmatrix} U & 0 \\ 0
& I_{\cU} \end{bmatrix}.
$$
\begin{theorem}
\label{C:gencanuni}
Any  observable weakly coisometric realization of a Schur
function $S\in \cS_{d}(\cU, \cY)$ is unitarily equivalent to
some generalized functional-model realization of $S$.
              \end{theorem}
\begin{proof} Let
$S(\blam)=D+\widetilde{C}(I_{\cX}-Z(\blam)\widetilde{A})^{-1}
Z(\blam)\widetilde{B}$ be an observable weakly coisometric realization of
$S$ with the state space $\widetilde \cX$. Then
$\cH(K_S)=\cH(K_{\widetilde{C},\widetilde{\bA}})$ by Proposition
\ref{R:1.3}. The observability operator  ${\mathcal
O}_{\widetilde{C}, \widetilde{\bA}}: \, x\to
\widetilde{C}(I_{\cX}-Z(\blam)\widetilde{A})^{-1}x$
associated with the contractive observable pair
$(\widetilde{C}, \widetilde{\bA})$ is isometric as an operator from
$\widetilde \cX$ into
$\cH(K_{\widetilde{C},\widetilde{\bA}})=\cH(K_S)$ by part (2)
in Theorem \ref{T:3-1.2nc}. Let us define the operator tuple ${\mathbf A}
=(A_{1}, \dots, A_{d})$ on the functional-model state space $\cX : =\cH(K_S)=
{\rm Ran} \, {\mathcal O}_{\widetilde{C}, \widetilde{\bA}}$ by
\begin{equation} \label{ups}
A_j {\mathcal O}_{\widetilde{C}, \widetilde{\bA}}=
{\mathcal O}_{\widetilde{C}, \widetilde{\bA}}\widetilde{A}_jx
\quad \text{for}\quad j=1,\ldots,d.
\end{equation}
Then for the generic element
$f(\blam)=\widetilde{C}(I_{\cX}-Z(\blam)\widetilde{A})^{-1}x$ of
$\cH(K_S)$ we have
\begin{eqnarray*}
f(\blam)-f(0)&=&\widetilde{C}(I-Z(\blam)\widetilde{A})^{-1}x
-\widetilde{C}x \\
&=&\widetilde{C}(I-Z(\blam)\widetilde{A})^{-1}Z(\blam)\widetilde{A}x \\
          &=&\sum_{j=1}^d \lambda_j
\widetilde{C}(I-Z(\blam)\widetilde{A})^{-1}\widetilde{A}_jx \\
          &=&\sum_{j=1}^d  \lambda_j \cdot({\cO}_{\widetilde{C},
          \widetilde{\bA}}\widetilde{A}_jx)(\blam)\\
          &=&\sum_{j=1}^d  \lambda_j\cdot
          (A_j{\cO}_{\widetilde{C}, \widetilde{\bA}}x)(\blam)
          =\sum_{j=1}^d  \lambda_j\cdot (A_j f)(\blam)
          \end{eqnarray*}
which means that the operators $A_1,\ldots,A_d$ solve the Gleason problem
on $\cH(K_S)$. For the same generic element $f(\blam)$ of $\cH(K_S)$
and for the operator $C: \, \cH(K_S)\to \cY$ defined as in
\eqref{dop8} we  also have
$$
C{\mathcal O}_{\widetilde{C}, \widetilde{\bA}}x=Cf=f(0)=\widetilde{C}x
$$
and, since the vector $x\in\cX$ is arbitrary, it follows that
\begin{equation}\label{upsa}
C{\mathcal O}_{\widetilde{C}, \widetilde{\bA}}=\widetilde{C}.
\end{equation}
Now we let
\begin{equation}\label{upsb}
B_j:={\mathcal O}_{\widetilde{C}, \widetilde{\bA}}
\widetilde{B}_j\quad\text{for} \quad j=1,\ldots,d.
\end{equation}
It is readily seen that $B_j$ maps $\cU$ into $\cH(K_S)$ and it follows
from \eqref{ups}, \eqref{upsa} and \eqref{upsb} that the realization
$\bU=\sbm{A & B \\ C & D}$ is unitarily
equivalent to the original realization
$\widetilde{\bU}=\sbm{\widetilde{A} & \widetilde{B} \\
\widetilde{C} & D}$ via the unitary operator $\cO_{\widetilde C,
\widetilde \bA} \colon \cX \to \cH(K_{S})$:
$$
\begin{bmatrix} A & B \\ C & D \end{bmatrix} \begin{bmatrix}
      \cO_{\widetilde C, \widetilde \bA} & 0 \\ 0 & I_{\cU} \end{bmatrix}
      = \begin{bmatrix} \oplus_{k=1}^{d} \cO_{\widetilde C, \widetilde
      \bA} & 0 \\ 0 & I_{\cY }\end{bmatrix} \begin{bmatrix} \widetilde
      A & \widetilde B \\ \widetilde C & D \end{bmatrix}.
$$
Therefore this realization $\bU$ is also weakly
coisometric. Also it is a generalized functional-model realization
since the state space $\cX$ is the functional-model state space
$\cH(K_{S})$, the output operator $C$ is given by evaluation at $0$, and
the state space operators
$A_1,\ldots, A_d$ on $\cX = \cH(K_{S})$ solve the Gleason problem in
$\cH(K_{S})$.
\end{proof}

As we have already seen, a Schur class function $S\in\cS_d(\cU,\cY)$
can admit more than one (not unitarily equivalent) weakly
coisometric realizations of the form \eqref{1.5a} with the same
$A_1,\ldots,A_d$ and $C$. Theorem \ref{T:noncomreal} indicates another
source for nonuniqueness: the kernel $K_S$ can be represented in the
form $K_{C,\bA}$ in more than one way, or equivalently, the Gleason
problem for the space $\cH(K_S)$ may have contractive solutions that
are not unitarily equivalent. A description of all contractive solutions
of the Gleason problem lies beyond the scope of this paper and will be
presented elsewhere. Here we present an example showing that the
nonuniqueness of the representing pair $(C, \bA)$ indeed may occur.
\begin{example}
\label{E:3.3}
{\rm Let us introduce the matrices
\begin{equation}\label{aug1}
C=\begin{bmatrix}\frac{1}{2} & 0 & 0 \end{bmatrix},  \quad
A_{0,1}=\begin{bmatrix} 0 & \frac{1}{4} & 0 \\ 0 & 0 & 0 \\ \frac{1}{2}
&0 &0\end{bmatrix},  \quad
A_{0,2}=\begin{bmatrix} 0 & 0 & \frac{1}{4} \\ \frac{1}{2} & 0 & 0
       \\ 0 & 0 & 0\end{bmatrix},
\end{equation}
\begin{eqnarray}
B_{0,1}&=&\begin{bmatrix}0 & \frac{\sqrt{15}}{4} & 0 & 0 & 0 & 0 & 0 \\
0 & 0 & 1& 0 & 0 & 0 & 0 \\
-\frac{1}{2\sqrt{3}} & 0 & 0 & \sqrt{\frac{2}{3}} & 0 & 0 & 0
\end{bmatrix},
\label{aug2}\\
B_{0,2}&=&\begin{bmatrix}0 & 0 & 0 & 0 & \frac{\sqrt{15}}{4} & 0 & 0 \\
-\frac{1}{2\sqrt{3}} &  0 & 0 & -\frac{1}{\sqrt{6}} & 0 &
\frac{1}{\sqrt{2}} & 0 \\
0 & 0 & 0 & 0 & 0 & 0 & 1 \end{bmatrix},
\label{aug3}\\
D&=&\begin{bmatrix} \frac{\sqrt{3}}{2} & 0 & 0 & 0 & 0 & 0 &
0\end{bmatrix}\label{aug4}
\end{eqnarray}
so that the $7\times 10$ matrix
$$
{\bf U}_0=\begin{bmatrix}A_{0,1} & B_{0,1} \\ A_{0,2} & B_{0,2} \\
C & D\end{bmatrix}
$$
is coisometric. Then the characteristic function
\begin{equation}\label{aug6}
S(\blam)
=D+C(I-\lambda_1A_{0,1}-\lambda_2A_{0,2})^{-1}
(\lambda_1B_{0,1}+\lambda_2B_{0,2})
\end{equation}
of the colligation ${\bf U}_0$
belongs to the Schur class $\cS_2(\C^7,\C)$. It is readily seen that
\begin{equation}\label{aug7}
C(I-\lambda_1A_{0,1}-\lambda_2A_{0,2})^{-1}=
\frac{1}{2} \begin{bmatrix}
\frac{4}{4-\lambda_1\lambda_2} & \frac{\lambda_1}{4-\lambda_1\lambda_2} &
\frac{\lambda_2}{4-\lambda_1\lambda_2}\end{bmatrix}
\end{equation}
which being substituted along with \eqref{aug3}--\eqref{aug4} into
\eqref{aug6} gives the explicit formula
\begin{equation}\label{aug8}
S(\blam)=\frac{1}{2(4-\lambda_1\lambda_2)}
\begin{bmatrix}\frac{12-4\lambda_1\lambda_2}{\sqrt{3}} &
\sqrt{15}\lambda_1 & \lambda_1^2 & \frac{\lambda_1\lambda_2}{\sqrt{6}}
& \sqrt{15}\lambda_2 & \frac{\lambda_1\lambda_2}{\sqrt{2}} &
\lambda_2^2\end{bmatrix}.
\end{equation}
By \eqref{aug7}, identity
$$
C(I-\lambda_1A_{0,1}-\lambda_2 A_{0,2})^{-1}\begin{bmatrix} x_1 \\ x_2
\\x_3\end{bmatrix}=\frac{1}{2} \cdot
\frac{4x_1+x_2\lambda_1+x_3\lambda_2}{4-\lambda_1\lambda_2}
\equiv 0
$$
implies $x_1=x_2=x_3=0$ and therefore the pair $(C,\bA_0)$ is observable.
Thus, representation \eqref{aug6} is a coisometric (and
therefore, also weakly coisometric) observable realization of
the function $S\in\cS_2(\C^7,\C)$ given by \eqref{aug8}. Then we also
have
\begin{eqnarray}
K_S(\blam,\bzeta)&=&C(I-\lambda_1A_{0,1}-\lambda_2A_{0,2})^{-1}
(I-\bar{\zeta}_1A_{0,1}^*-\bar{\zeta}_2A_{0,2}^*)^{-1}C^*\nonumber\\
&=&K_{C,\bA_0}(\blam,\bzeta).\label{aug9}
\end{eqnarray}
Now let us consider the matrices
\begin{equation}\label{aug10}
A_{\gamma,1}=\begin{bmatrix} 0 & \frac{1}{4} & 0 \\ 0 & 0 & 0 \\
\frac{1}{2}+\gamma
&0 &0\end{bmatrix}\quad\mbox{and}\quad
A_{\gamma,2}=\begin{bmatrix} 0 & 0 & \frac{1}{4} \\ \frac{1}{2}-\gamma & 0
& 0     \\ 0 & 0 & 0\end{bmatrix}
\end{equation}
where $\gamma\in\C$ is a parameter, and note that
$$
C(I-\lambda_1A_{\gamma,1}-\lambda_2A_{\gamma,2})^{-1}=
\frac{1}{2} \begin{bmatrix}
\frac{4}{4-\lambda_1\lambda_2} & \frac{\lambda_1}{4-\lambda_1\lambda_2} &
\frac{\lambda_2}{4-\lambda_1\lambda_2}\end{bmatrix}
$$
for every $\gamma$. In particular, the pair  $(C,\bA_\gamma)$ is
observable for every $\gamma$. The latter equality together with
\eqref{aug9} gives
\begin{equation}\label{aug11}
K_S(\blam,\bzeta)=K_{C,\bA_\gamma}(\blam,\bzeta).
\end{equation}
Now pick any $\gamma$ so that $|\gamma|<\frac{1}{2\sqrt{2}}$. As it is
easily seen, the latter inequality is equivalent to the pair $(C,{\bf
A}_\gamma)$ being contractive. Thus, we have a Schur class function $S$
and  a contractive pair $(C,{\bf A}_\gamma)$ such that equality
\eqref{aug11}
holds. Then by Theorem \ref{T:CAtoS}, there exist operators $B_{\gamma,1}$
and $B_{\gamma,2}$ so that the operator
$$
{\bf U}_\gamma=\begin{bmatrix}A_{\gamma,1} & B_{\gamma,1} \\ A_{\gamma,2}
& B_{\gamma,2} \\ C & D\end{bmatrix}
$$
is weakly coisometric and $S$ can be realized as
$$
S(\blam)
=D+C(I-\lambda_1A_{\gamma,1}-\lambda_2A_{\gamma,2})^{-1}
(\lambda_1B_{\gamma,1}+\lambda_2B_{\gamma,2}).
$$
It remains to note that the pairs $(C,\bA_{\gamma})$ and
$(C,\bA_{\gamma'})$ are not unitarily equivalent (which is shown by
another elementary calculation) unless $\gamma=\gamma'$.}
\end{example}

\section{Overlapping spaces}
\label{S:over}

          The subspaces ${\mathcal D}^{\perp}$ and ${\mathcal
          R}_{V}^{\perp}$ as described in
          \eqref{canonical3}, \eqref{canonical2} are
          particular examples of a general notion of {\em overlapping
          spaces} appearing in the theory of reproducing kernel Hilbert
          spaces as developed by de Branges and Rovnyak \cite{dbr1, dbr2}.
          In general, suppose that ${\mathbb M} = {\mathbb M}(\blam,
          \bzeta)$ is a positive kernel on $\Omega \times \Omega$
          with values in ${\mathcal L}(\cX)$ (for some Hilbert space $\cX$)
          inducing a reproducing kernel Hilbert space $\cH({\mathbb M})$ of
          $\cX$-valued functions via the Aronszajn construction, and
          suppose that $F$ is a function on $\Omega$ with values equal to
          operators from $\cX$ to another Hilbert space $\cX'$.
          (In our application, of course, we will take $\Omega =
{\mathbb B}^{d}$.)
          Then
          \begin{equation} \label{MF}
         {\mathbb M}_{F}(\blam, \bzeta): = F(\blam) {\mathbb
          M}(\blam, \bzeta) F(\bzeta)^{*}
          \end{equation}
          is also a positive kernel on
          $\Omega \times \Omega$ with values in ${\mathcal L}(\cX')$
          inducing a reproducing kernel Hilbert space $\cH({\mathbb
          M}_{F})$ of $\cX'$-valued functions on $\Omega$.
          The sets of finite linear combinations of kernel functions
          \begin{align*}
         & {\mathcal S}_{{\mathbb M}}: = \left\{ \sum_{k=1}^{N} {\mathbb
          M}(\cdot, \bzeta_{k}) x_{k} \colon \bzeta_{k} \in \Omega,\, x_{k}
          \in \cX, \, N = 1,2,3, \cdots \right\}, \\
         & {\mathcal S}_{{\mathbb M}_{F}}: = \left \{ \sum_{k=1}^{N} {\mathbb
          M}_{F}(\cdot, \bzeta_{k}) x_{k} \colon \bzeta_{k} \in \Omega,\, x_{k}
          \in \cX, \, N = 1,2,3, \cdots \right\}
         \end{align*}
         form dense sets in $\cH({\mathbb M})$ and $\cH({\mathbb M}_{F})$
         respectively.  Moreover, the computation
         \begin{align*}
             \left\| \sum_{k=1}^{N} {\mathbb M}_{F}(\cdot, \bzeta_{k}) x'_{k}
             \right\|^{2} & =
             \sum_{k,\ell=1}^{N} \left\langle {\mathbb M}_{F}(\cdot, \bzeta_{k})
             x_{k}', {\mathbb M}_{F}(\cdot, \bzeta_{\ell}) x'_{\ell}
             \right\rangle_{\cH({\mathbb M}_{F})} \\
             & = \sum_{k,\ell=1}^{N} \langle {\mathbb M}_{F}(\bzeta_{\ell},
             \bzeta_{k}) x'_{k}, x'_{\ell} \rangle_{\cX'} \\
             & = \sum_{k,\ell=1}^{N} \langle {\mathbb M}(\bzeta_{\ell},
             \bzeta_{k}) F(\bzeta_{k})^{*}x'_{k}, F(\bzeta_{\ell})^{*}
             x'_{\ell} \rangle_{\cX} \\
         & = \sum_{k,\ell=1}^{N} \langle {\mathbb M}(\cdot, \bzeta_{k})
             F(\bzeta_{k})^{*} x_{k}', {\mathbb M}(\cdot, \bzeta_{\ell})
             F(\bzeta_{\ell})^{*} x'_{\ell} \rangle_{\cH({\mathbb M})} \\
             & = \left\| \sum_{k=1}^{N} {\mathbb M}(\cdot, \bzeta_{k})
             F(\bzeta_{k})^{*} x'_{k} \right\|^{2}_{\cH({\mathbb M})}
         \end{align*}
         shows that the map
         $$ \Psi \colon {\mathbb M}_{F}(\cdot, \bzeta) x' \mapsto {\mathbb
         M}(\cdot, \bzeta) F(\bzeta)^{*} x'
         $$
         (for $\bzeta \in \Omega$ and $x' \in \cX'$) extends by linearity
         and continuity to define an isometry, still called $\Psi$, from
         $\cH({\mathbb M}_{F})$ into $\cH({\mathbb M})$.  Another
         computation, where $f \in \cH({\mathbb M})$, $\bzeta \in \Omega$
         and $x' \in \cX'$,
         \begin{align*}
             \langle \Psi^{*} f, {\mathbb M}_{F}(\cdot, \bzeta) x' \rangle
             & = \langle f, \Psi {\mathbb M}_{F}(\cdot, \bzeta) x'
             \rangle_{\cH({\mathbb M})} \\
             & = \langle f, {\mathbb M}(\cdot, \bzeta) F(\bzeta)^{*} x'
             \rangle_{\cH({\mathbb M})} \\
             & = \langle f(\bzeta), F(\bzeta)^{*} x' \rangle_{\cX} \\
             & = \langle F(\bzeta) f(\bzeta), x' \rangle_{\cX'}
         \end{align*}
         shows that the adjoint of $\Psi$ is the multiplication operator
          $$\Psi^{*} =  M_{F} \colon f(\blam) \mapsto F(\blam) f(\blam).
          $$
          Since we saw above that $\Psi$ is an isometry, we conclude 
that $M_{F}$
          is a coisometry from $\cH({\mathbb M})$ onto $\cH({\mathbb
          M}_{F})$ and that $\cH({\mathbb M}_{F})$ can be
          characterized as
          $$ \cH({\mathbb M}_{F}) = \{ F \cdot f \colon f \in \cH({\mathbb
          M}) \}
          $$
          with norm given by
          \begin{align}
          \| F \cdot f \|_{\cH({\mathbb M}_{F})} & = \inf \{ \|f'
          \|_{\cH({\mathbb M})} \colon F(\blam) f'(\blam) = F(\blam)
          f(\blam) \text{ for all } \blam \in \Omega \} \notag \\
          & =  \| Q f \|_{\cH({\mathbb M})}\quad \text{ where}\quad Q =
          P_{(\operatorname{Ker}\, M_{F})^{\perp}}.
          \label{HMSnorm}
          \end{align}
          The associated {\em overlapping space} ${\mathcal
          L}(F, {\mathbb M})$ is defined to be
          $$
           {\mathcal L}(F, {\mathbb M}) = \operatorname{Ker}\, M_{F} \subset
           \cH({\mathbb M})
          $$
          with norm inherited from ${\mathcal H}({\mathbb M})$.  We then
          have the unitary identification map
          $$
          \Gamma := \begin{bmatrix} M_{F} \\ P_{\operatorname{Ker} \,M_{F}}
          \end{bmatrix}  \colon \; \cH({\mathbb M}) \to \begin{bmatrix}
\cH({\mathbb M}_{F}) \\ {\mathcal L}(F, {\mathbb M}) \end{bmatrix}.
          $$
          When there are canonical operators on $\cH({\mathbb M})$, it is
          often of interest to work out the induced canonical operators on
          $\cH({\mathbb M}_{F}) \oplus {\mathcal L}(F, {\mathbb M})$.  We
          discuss two particular instances here related to Proposition
          \ref{P:canonical-guys}; in these examples, $\Omega = {\mathbb
          B}^{d}$.
\begin{example}\label{E:1}
{\rm   Take ${\mathbb M}(\blam, \bzeta) = K_{S}(\blam, \bzeta)
        \otimes I_{{\mathbb C}^{d}}$ and $F(\blam) = Z(\blam)$.
Then $\cH({\mathbb M}) = \cH(K_{S})^{d}$ and the
             associated kernel ${\mathbb M}_{F}(\blam, \bzeta)$ is given by
             $$
             (K_{S}\otimes I_{{\mathbb C}^{d}})_{_Z} =
             \lam_{1}K_{S}(\blam, \bzeta) \overline{\zeta_{1}} + \cdots +
             \lam_{d}K_S(\blam, \bzeta) \overline{\zeta_{d}}
             $$
             with associated overlapping space ${\mathcal L}(F, {\mathbb
             M})$ given by
             $$
             {\mathcal L}(Z, K_{S}\otimes I_{{\mathbb C}^{d}}) =
             \{ f  \in \cH(K_{S})^{d} \colon Z(\blam) f(\blam) \equiv 0 \}
          $$
          Then ${\mathcal L}(Z, K_{S}\otimes I_{{\mathbb C}^{d}})$ is
         exactly the subspace ${\mathcal D}^{\perp}$ in \eqref{canonical2}.
          Thus $M_{Z} \colon f(\blam) \mapsto Z(\blam)f(\blam)$ is unitary
          from ${\mathcal D}$ onto $\cH((K_{S} \otimes I_{{\mathbb
C}^{d}})_{_Z})$.}
\end{example}
\begin{example}\label{E:2}
{\rm  Take
${\mathbb M}(\blam, \bzeta) = \begin{bmatrix} K_{S}(\blam, \bzeta)
          & 0 \\ 0 & I_{\cU} \end{bmatrix}$ and  $F(\blam) =
          \begin{bmatrix} I_{\cY} & S(\blam) \end{bmatrix}$.
         Then the associated kernel ${\mathbb M}_{F}(\blam, \bzeta)$ is
         $$
         \left(K_{S} \oplus I_{\cU}\right)_{\sbm{ I & S}}(\blam, \bzeta)
         = K_{S}(\blam, \bzeta) + S(\blam) S(\bzeta)^{*}
$$
while the  associated overlapping space ${\mathcal L}(F, {\mathbb M})$
is given by
$${\mathcal L}\left(\begin{bmatrix} I_{\cY} & S \end{bmatrix}, K_{S}
\oplus
I_{\cU}\right) = \left\{ \begin{bmatrix} h \\ u \end{bmatrix} \in
\begin{bmatrix} \cH(K_{S}) \\ \cU \end{bmatrix} \colon h(\blam) +
          S(\blam) u \equiv 0 \right\}
        $$
        and is exactly equal to the space ${\mathcal R}_{V}^{\perp}$ in
        \eqref{canonical2}.  Note that the space $\cU^{0}$ defined in
        \eqref{KerT22} is related to ${\mathcal L}(\begin{bmatrix}
        I_{\cY} & S \end{bmatrix}, K_{S}
        \oplus I_{\cU})$ according to
       $$
        \begin{bmatrix} 0 \\ \cU_{S}^{0}\end{bmatrix} =
        {\mathcal L}(\begin{bmatrix} I_{\cY} & S \end{bmatrix}, K_{S} \oplus
I_{\cU})    \bigcap \begin{bmatrix} 0 \\ \cU \end{bmatrix}.
        $$
}\end{example}

          Overlapping spaces are usually
	considered only for the case where $F$ and ${\mathbb M}$ have
	the special form
	$$ F(\blam, \bzeta) = \begin{bmatrix} F_{1}(\blam) &
	F_{2}(\blam) \end{bmatrix}, \qquad {\mathbb M}(\blam, \bzeta)
	= \begin{bmatrix} {\mathbb M}_{1}(\blam, \bzeta) & 0 \\ 0 &
	{\mathbb M}_{2}(\blam, \bzeta) \end{bmatrix}
	$$
	(see \cite{dbr1, dbr2}), but the case of any finite number
	(or even a continuum) of such positive kernels ${\mathbb
	M}_{s}(\blam, \bzeta)$ has come up in
	some applications (see \cite{Cuntz-rep}).

\end{document}